\documentclass[a4paper]{article}
\usepackage{amsmath,amssymb,amscd,amsthm}
\title{Representation of the Gauss hypergeometric function by multiple polylogarithms and relations of multiple zeta values}
\author{Shu OI}

\DeclareFontEncoding{OT2}{}{}
\DeclareFontSubstitution{OT2}{wncyr}{m}{n}
\DeclareSymbolFont{cyss}{OT2}{wncyss}{m}{n}
\DeclareSymbolFont{cyr}{OT2}{wncyr}{m}{n}
\DeclareMathSymbol{\sh}{\mathbin}{cyss}{`x}

\newcommand{\Li}{\operatorname{Li}}
\newcommand{\vect}[2]{{\begin{pmatrix}{#1}\\{#2}\end{pmatrix}}}
\newcommand{\reg}{\operatorname{reg}}
\newcommand{\dist}{\operatorname{dist}}

\newcommand{\R}{{\mathbf R}}
\newcommand{\Q}{{\mathbf Q}}
\newcommand{\Z}{{\mathbf Z}}
\newcommand{\N}{{\mathbf N}}
\newcommand{\bP}{{\mathbf P}}

\newcommand{\bk}{{\mathbf k}}

\newcommand{\bw}{{\mathbf w}}
\newcommand{\bx}{\text{\boldmath $x$}}
\newcommand{\bmu}{\text{\boldmath $\mu$}}

\newcommand{\fH}{{\mathfrak H}}

\newcommand{\fU}{{\mathfrak U}}

\newcommand{\ds}{\displaystyle}

\theoremstyle{definition}
\newtheorem{thm}{Theorem}
\newtheorem{cor}[thm]{Corollary}
\newtheorem{prop}[thm]{Proposition}
\newtheorem{lem}[thm]{Lemma}

\begin{document}


\maketitle

\begin{abstract}
We describe a solution of the Gauss hypergeometric equation, $F(\alpha,\beta,\gamma;z)$ by power series in paramaters $\alpha,\beta,\gamma$ whose coefficients are $\Z$ linear combinations of multiple polylogarithms.
And using the representation and connection reletions of solutions of the hypergeometric equation, we show some relations of multiple zeta values.
\end{abstract}

\section{Introduction} \label{sec:Introduction}

In introduction of Aomoto-Kita\cite{AoK}, they say that the Gauss hypergeometric equation:
\begin{equation}
z(1-z)\frac{d^2w}{dz^2} + (\gamma - (\alpha+\beta+1)z) \frac{dw}{dz} - \alpha \beta w = 0 \label{HGeq}
\end{equation}
can be solved by iterated integral and the solutions can be described using multiple polylogarithms, but concrete representation is not appeared.

Multiple polylogarithms\;\;$\Li_{\bk}(z)$ ($\bk=(k_1,\ldots,k_n)$ is a sequence of natural numbers) are defined by
\begin{equation}
\Li_{\bk}(z) = \sum_{m_1>m_2>\cdots>m_n>0}\frac{z^{m_1}}{m_1^{k_1}\cdots m_n^{k_n}}. \label{MPLdef}
\end{equation}
which are holomorphic functions defined on $|z|<1$.

In particular if $k_1>1$, these functions converge as $z \to 1$,
\begin{equation}
\Li_{\bk}(1) = \zeta(k_1,\ldots,k_n) = \sum_{m_1>m_2>\cdots>m_n>0}\frac{1}{m_1^{k_1}\cdots m_n^{k_n}} \quad \text{(multiple zeta values).} \label{MZVdef}
\end{equation}

Connection relations of solutios of the hypergeometric equation are written using the gamma functions explicitly. Then if solutions of the hypergeometric equation are described by multiple polylogarithms, we will get some ralations of multiple zeta values by comparing connection relation and limits of the solutios as $z$ tend to each singular point.

On this problem, Ohno-Zagier\cite{OZ} show that a generating function of sum of multiple polylogarithms fixed "weight", "depth" and "height" satisfies the hypergeometric equation and show relations of multipe zeta values using Gauss formula about special values of the hypergeometric function: a solution of the hypergeometric equation regular on $z=0$. 

\paragraph{}

In this paper, we try to solve the hypergeometric equation by iterated integral and express solutions by multiple polylogarithms concretely. And then comparing the representations and connection relation of solutions of hypergeometric equation, we show some ralations of multiple polylogarithms and multiple zeta values.

The results include the formula shown by Ohno-Zagier\cite{OZ} and partially the duality formula. The results give natural interpretation of Ohno-Zagier formula by method of solving the equation directly.

\paragraph{}

Regarding the hypergeometric equation as easiest case of Knizhnik-Zamolodchikov(KZ) equation, we can consider the method of representing solutions by multiple polylogarithms is applied to more general KZ equations. A position of this paper is the first step of these application.

\section{Preliminary} \label{sec:Preliminary}

\subsection{The Gauss hypergeometric equation and its solutions} \label{subsec:Pre:HGE}

The Gauss hypergeometric equation is a Fuchsian differential equation of second order which have three reguler singuler points $0,1,\infty$ in $\bP^1$. The equation is defined by 
\begin{equation}
z(1-z)\frac{d^2w}{dz^2} + (\gamma - (\alpha+\beta+1)z) \frac{dw}{dz} - \alpha \beta w = 0.
\end{equation}
Now $\alpha,\beta,\gamma$ are parameters. Each solution of the equation can be connected on $\bP-\{0,1,\infty\}$ as many-valued analytic function.

The equation has a series solution which is regular at the origin:
\begin{equation}
F(\alpha,\beta,\gamma;z) = \sum_{n=0}^{\infty} \frac{(\alpha)_n(\beta)_n}{(\gamma)_n n!} z^n.
\end{equation}
On condition that $(\alpha)_n=\alpha(\alpha+1)(\alpha+2)\cdots(\alpha+n-1)$.

This series converges on $|z|<1$ if $\gamma \not\in \Z_{\le 0}$. The function is named the Gauss hypergeometric function.

All solutions of the hypergeometric equation \eqref{HGeq} in neighborhood of $0,1$ and $\infty$ can be written by the hypergeometric function.

Two linearly independental solutions in neighborhood of $0$ are
\begin{align}
\varphi^{(0)}_0(z) &= F(\alpha,\beta,\gamma;z)\\
\intertext{and}
\varphi^{(0)}_1(z) &= z^{1-\gamma}F(\alpha-\gamma+1,\beta-\gamma+1,2-\gamma;z).
\end{align}

Solutions in neighborhood of $1$ are
\begin{align}
\varphi^{(1)}_0(1-z) &= F(\alpha,\beta,\alpha+\beta-\gamma+1;1-z)\\
\intertext{and}
\varphi^{(1)}_1(1-z) &= (1-z)^{\gamma-\alpha-\beta}F(\gamma-\alpha,\gamma-\beta,\gamma-\alpha-\beta+1;1-z).
\end{align}

Solutions in neighborhood of $\infty$ are
\begin{align}
\varphi^{(\infty)}_0(1/z) &= z^{-\alpha}F(\alpha,\alpha-\gamma+1,\alpha-\beta+1;1/z)\\
\intertext{and}
\varphi^{(\infty)}_1(1/z) &= z^{-\beta}F(\beta,\beta-\gamma+1,\beta-\alpha+1;1/z).
\end{align}

These solutions are related by analytic continuation. The connection coefficients are given by the following fomula.

We set fundamental solution matrix on $i=0,1,\infty$ by
\begin{equation}
\Phi_i=\begin{pmatrix}\varphi^{(i)}_0&\varphi^{(i)}_1\\\frac{1}{\beta}z\frac{d}{dz}\varphi^{(i)}_0&\frac{1}{\beta}z\frac{d}{dz}\varphi^{(i)}_1\end{pmatrix},
\end{equation}
the connection relations are written as
\begin{align}
\Phi_0 & = \Phi_1 C^{01}, \;\;C^{01}=\begin{pmatrix} \frac{\ds \Gamma(\gamma)\Gamma(\gamma-\alpha-\beta)}{\ds \Gamma(\gamma-\alpha)\Gamma(\gamma-\beta)} & \frac{\ds \Gamma(2-\gamma)\Gamma(\gamma-\alpha-\beta)}{\ds \Gamma(1-\alpha)\Gamma(1-\beta)} \\ \frac{\ds \Gamma(\gamma)\Gamma(\alpha+\beta-\gamma)}{\ds \Gamma(\alpha)\Gamma(\beta)} & \frac{\ds \Gamma(2-\gamma)\Gamma(\alpha+\beta-\gamma)}{\ds \Gamma(\alpha-\gamma+1)\Gamma(\beta-\gamma+1)} \end{pmatrix} \label{connection01}\\
\Phi_0 & = \Phi_{\infty} C^{0\infty}, \;\;C^{0\infty}=\begin{pmatrix} e^{-\pi i \alpha} \frac{\ds \Gamma(\gamma)\Gamma(\beta-\alpha)}{\ds \Gamma(\beta)\Gamma(\gamma-\alpha)} & e^{\pi i (\gamma-\alpha-1)} \frac{\ds \Gamma(2-\gamma)\Gamma(\beta-\alpha)}{\ds \Gamma(\beta-\gamma+1)\Gamma(1-\alpha)} \\ e^{-\pi i \beta} \frac{\ds \Gamma(\gamma)\Gamma(\alpha-\beta)}{\ds \Gamma(\alpha)\Gamma(\gamma-\beta)} & e^{\pi i (\gamma-\beta-1)} \frac{\ds \Gamma(2-\gamma)\Gamma(\alpha-\beta)}{\ds \Gamma(\alpha-\gamma+1)\Gamma(1-\beta)} \end{pmatrix} \label{connection0infty}
\end{align}

The $(1,1)$ component of \eqref{connection01} tends to Gauss's formula as $z \to 1$:
\begin{equation}
F(\alpha,\beta,\gamma;1) = \frac{\ds \Gamma(\gamma)\Gamma(\gamma-\alpha-\beta)}{\ds \Gamma(\gamma-\alpha)\Gamma(\gamma-\beta)} \qquad (\text{when } \Re(\gamma-\alpha-\beta)>0), \label{GAUSSformula}
\end{equation}

\subsection{Multiple zeta values and multiple polylogarithms} \label{subsec:Pre:MZV}

Multiple zeta values are real numbers defined by 
\begin{equation}
\zeta(k_1,\ldots,k_n) = \sum_{m_1>m_2>\cdots>m_n>0}\frac{1}{m_1^{k_1}\cdots m_n^{k_n}}.
\end{equation}
Now $k_1,\ldots,k_n$ are natural numbers and $k_1 \ge 2$. The range of summation is all integers $m_1>m_2>\cdots>m_n>0$. On the condition this serieses converge to real numbers. If $n=1$, $\zeta(k_1)$ are Riemann zeta values.

Multiple polylogarithms are defined by the following power serieses:
\begin{equation}
\Li_{\bk}(z) = \sum_{m_1>m_2>\cdots>m_n>0}\frac{z^{m_1}}{m_1^{k_1}\cdots m_n^{k_n}}.
\end{equation}
Now $\bk=(k_1,\ldots,k_n)$ is a sequence of natural numbers. The convergence radius of these serieses are 1. If $k_1 \ge 2$, these serieses converge to multiple zeta values $\zeta(k_1,\ldots,k_n)$ as $z \to 1$.

On multiple zeta values and multiple polylogarithms, $k_1+\cdots+k_n$, $n$ and $\#\{i|k_i \ge 2\}$ are named {\it weight},{\it depth} and {\it height} respectively.
	
$\Li_{\bk}(z)$ satisfy differential relations as follows:
\begin{equation}
\frac{d \Li_{k_1,\ldots,k_n}(z)}{dz} = 
	\begin{cases}
	\frac{1}{1-z}\Li_{k_2,\ldots,k_n}(z) & \text{if $k_1=1$}\\
	\frac{1}{z}\Li_{k_1-1,\ldots,k_n}(z) & \text{if $k_1>1$}
	\end{cases}
\end{equation}
and $\Li_{\bk}(z)$ have iterated integral representations:
\begin{equation}
\Li_{k_1,\ldots,k_n}(z) = \int_0^z \underbrace{\frac{dt}{t} \circ \frac{dt}{t} \circ \cdots \circ \frac{dt}{t}}_{k_1-1 \text{ times}} \circ \frac{dt}{1-t} \circ \underbrace{\frac{dt}{t} \circ \frac{dt}{t} \circ \cdots \circ \frac{dt}{t}}_{k_2-1 \text{ times}} \circ \frac{dt}{1-t} \circ \cdots \circ \underbrace{\frac{dt}{t} \circ \frac{dt}{t} \circ \cdots \circ \frac{dt}{t}}_{k_n-1 \text{ times}} \circ \frac{dt}{1-t}. \label{rep_ite_MPL}
\end{equation}

Now notation $\ds \int_0^z \omega_1(t) \circ \omega_2(t) \circ \cdots \circ \omega_n(t)$ (each $\omega_i(t)$ is a differntial form of $t$) means\\
$\ds \int_0^z \omega_1(t_1) \int_0^{t_1} \omega_2(t_2) \cdots \int_0^{t_{n-1}}\omega_n(t_n)$.

For example, if $n=1$ and $k_1=1$,
\begin{equation*}
\Li_1(z)=\sum_{m=1}^{\infty}\frac{z^m}{m} = \int_0^z \frac{dt}{1-t}= -\log(1-z)
\end{equation*}

Multiple polylogarithms $\Li_{\bk}(z)$ can be connected analytically on $\bP^1-\{0,1,\infty\}$ by integral representation \eqref{rep_ite_MPL} and define many-valued analytic functions.

By changing variables on iterated integral representation of multiple zeta values(case of $z=1$ as \eqref{rep_ite_MPL}), we get {\it the duality formula of multiple zeta values}:
\begin{multline}
\zeta(a_1+1,\underbrace{1,\ldots,1}_{b_1-1 \text{ times}},a_2+1,\underbrace{1,\ldots,1}_{b_2-1 \text{ times}},\cdots,a_r+1,\underbrace{1,\ldots,1}_{b_r-1 \text{ times}}) \\
= \zeta(b_r+1,\underbrace{1,\ldots,1}_{a_r-1 \text{ times}},b_{r-1}+1,\underbrace{1,\ldots,1}_{a_{r-1}-1 \text{ times}},\cdots,b_1+1,\underbrace{1,\ldots,1}_{a_1-1 \text{ times}}). \label{MZVduality}
\end{multline}

\subsection{Ring of noncommutative polynomials $\Q\langle x,y\rangle$ and shuffle product} \label{subsec:Pre:Q<xy>}

We consider a noncommutative polynomial algebra $\fH =\Q\langle x,y\rangle$ in two indeterminant $x,y$ and subalgebras $\fH^1=\Q + \fH y$ and $\fH^0=\Q + x\fH y$.

We define a multiple polylogarithm with respect to monomial $x^{k_1-1}y \cdots x^{k_n-1}y \in \fH^1$ as
\begin{equation}
\Li(x^{k_1-1}y \cdots x^{k_n-1}y;z)=\Li_{k_1,\ldots,k_n}(z)
\end{equation}
and to each $\bw \in \fH^1$, $\Li(\bw;z)$ defined by above and $\Q$ linearity.

If $\bw \in \fH^1$, we get
\begin{align}
\frac{d}{dz}\Li(x\bw;z) &= \frac{1}{z}\Li(\bw;z) \notag\\
\frac{d}{dz}\Li(y\bw;z) &= \frac{1}{1-z}\Li(\bw;z). \label{MPLdiffeq}
\end{align}

Especially if $k_1 \ge 2$, $\Li(x^{k_1-1}y \cdots x^{k_n-1}y;1)=\zeta(k_1,\ldots,k_n)$.

We define anti-automorphism $\tau:\fH \to \fH$ as $\tau(x)=y,\tau(y)=x$. Using this notation, we can denote duality formula \eqref{MZVduality} as
\begin{equation}
\Li(\bw;1)=\Li(\tau(\bw);1)
\end{equation}
on $\bw \in \fH^0$.

\paragraph{}

Shuffle product $\sh$ on $\fH$ is defined by
\begin{align*}
&\sh \text{ is $\Q$-bilinear}\\
&1 \sh \bw = \bw \sh 1 = \bw \quad \bw \in \fH\\
&a_1\bw_1 \sh a_2\bw_2 = a_1(\bw_1 \sh a_2\bw_2) + a_2(a_1\bw_1 \sh \bw_2) \quad a_1,a_2 = x \text{ or } y,\;\; \bw_1,\bw_2 \in \fH
\end{align*}

According to Reutenauer\cite{R} and Ihara-Kaneko\cite{IK}, $(\fH,\sh)$ is a commutative algebra and $\fH^1, \fH^0$ are subalgebras. Moreover $\fH \cong \fH^1[x]=\oplus_{n=0}^{\infty} \fH^1 \sh x^{\sh n}$ and the map $\reg^1 : \fH \to \fH^1$ is defined by $\reg^1(\bw)=(\text{constant term of }\bw \in \fH^1[x])$.

The map $\reg^{1}$ satisfies 
\begin{align}
\bw x^n &= \sum_{j=0}^n \reg^{1}(\bw x^{n-j}) \sh x^j \quad \text{for } \bw \in \fH^{1}\\
\reg^1(\bw yx^n)&=(-1)^n(\bw \sh x^n)y \quad \text{for } n \ge 0, \bw \in \fH
\end{align}
shown by Ihara-Kaneko\cite{IK}.

By property of iterated integral, we get
\begin{equation}
\Li(\bw_1 \sh \bw_2;z)=\Li(\bw_1;z)\Li(\bw_2;z).
\end{equation}

\paragraph{}

We expand index of multiple polylogarithms from $\fH^{1}$ to $\fH$ by
\begin{equation}
\Li(\bw x^n;z) = \sum_{j=0}^n \Li(\reg^{1}(\bw x^{n-j});z)\frac{\log^j z}{j!}
\end{equation}
for $\bw \in \fH^{1}$ and $\Q$ linearity. This expantion is equal to expantion as $\Li(x;z)=\log z$ and $\sh$ homomorphism.

According to Okuda\cite{O}, this expantion satisfies the same differential equation \eqref{MPLdiffeq}.

We redefine weight, depth and height of $\bw \text{: monomial } \in \fH$ as
\begin{align*}
\text{weight of }\bw &= |\bw| = \text{number of letters in } \bw\\
\text{depth of }\bw &= d(\bw) = \text{number of letter '$y$' in } \bw\\
\text{height of }\bw &= h(\bw) = (\text{number of juxtaposition '$yx$' in } \bw) +1
\end{align*}
and weight,depth and height of multiple polylogarithms and multiple zeta values as same sence. This definition of height is not equal to previous definition if $\bw \not\in \fH^0$, but the definition is more natural in this paper.

Detail of this and previous subsection is written in Arakawa-Kaneko\cite{ArK}.

\subsection{Analytic property of multiple polylogarithms}\label{subsec:Pre:Anal_MPL}

According to Lappo-Danilevsky\cite{L}, we define $\Li_{p,C}(\bw;z)$ as
\begin{align}
\Li_{p,C}(x\bw;z)&=\int_{p,C}^z \frac{\Li_{p,C}(\bw;z)}{z}dz\\
\Li_{p,C}(y\bw;z)&=\int_{p,C}^z \frac{\Li_{p,C}(\bw;z)}{1-z}dz\\
\Li_{p,C}(1;z)&=1
\end{align}
and $\Q$ linearity by $\bw$. Now $\bw \in \fH,\;\; z, p \in \fU = \text{universal covering on } \bP-\{0,1,\infty\}$ and $C$ is a integral path from $p$ to $z$ $\subset \fU$.

In this notation, our multiple polylogarithms $\Li(\bw;z)$ is written as
\begin{equation}
\Li(\bw;z) = \lim_{\substack{p \to 0\\ \arg{p}=0}} \Li_{p,C}(\bw;z)
\end{equation}
on $z \in \fU$.

The following lemmas was shown by Lappo-Danilevsky\cite{L}.

\begin{lem}
\begin{equation}
|\Li_{p,C}(\bw;z)| < \frac{1}{|\bw|!}\left(\frac{\sigma}{\delta}\right)^{|\bw|}. \label{HLestimate}
\end{equation}
Now $\ds \delta = \inf_{\substack{z' \in C\\ a=0,1}}\dist(z',a)$, $\sigma = \text{length of }C$.
\end{lem}

\begin{lem}
Let $\bw = a_1a_2\cdots a_r$, $a_i = x \text{ or } y$, $q \in C$, $C = C'+C''$, $C'$ is path from $p$ to $q$, $C''$ is path from $q$ to $z$,
\begin{equation}
\Li_{p,C}(\bw;z)=\sum_{i=0}^{r} \Li_{q,C''}(a_1\cdots a_i;z)\Li_{p,C'}(a_{i+1}\cdots a_r;q). \label{HLconnect}
\end{equation}
Now $\Li_{p,C'}(a_1\cdots a_0;q)=\Li_{q,C''}(a_{r+1}\cdots a_r;z)=1$.
\end{lem}

\begin{lem}
If $0<z<\frac{1}{2}$
\begin{equation}
|\Li(\bw;z)|<1 \quad \forall \bw \in \fH^1
\end{equation}
\end{lem}

\begin{proof}
Let $0<z<\frac{1}{2}$. We prove the lemma by induction on $|\bw|$.

If $|\bw|=1 \Leftrightarrow \bw=y$, $|\Li(\bw;z)|=|-\log(1-z)|< \log 2 < 1$.

Assuming the lemma satisfies by $\bw = x^{k_1-1}y\cdots x^{k_r-1}y$,
\begin{align*}
|\Li(x\bw;z)|&=\sum_{m_1>\cdots>m_r>0} \frac{z^{m_1}}{m_1^{k_1+1}\cdots m_r^{k_r}} < \sum_{m_1>\cdots>m_r>0} \frac{z^{m_1}}{m_1^{k_1}\cdots m_r^{k_r}} = |\Li(\bw;z)|<1\\
|\Li(y\bw;z)|&\le \int_0^z \frac{|\Li(\bw;z)|}{1-z} dz < \int_0^z \frac{1}{1-z} dz < 2 \cdot \frac{1}{2} = 1
\end{align*}

\end{proof}

Using the lemmas, we can estimate absolute values of multiple polylogarithms on $z \in \fU$.

\begin{prop}
Let $K \subset \fU : $ compact subset. There exist a constant $M_K$ depending only $K$,
\begin{equation}
|\Li(\bw;z)|<M_K \quad \forall \bw \in \fH
\end{equation}
\end{prop}

\begin{proof}

Let $\varepsilon = \min\{\text{(distance from \{0,1\} to $K$)},\frac{1}{2}\}$

$l(z;\varepsilon) = \ds \inf_{\substack{\text{$C$ : path from $\varepsilon$ to $z$ in $\fU$}\\\text{(distance from \{0,1\} to $C$)}\ge \varepsilon}} (\text{length of }C)$

$\ds l(z;K)=\max_{z \in K}l(z;\varepsilon)$.

We assume $z \in K$, $\bw = a_1a_2 \cdots a_r \in \fH^1 \;\;(\Leftrightarrow a_r=y)$ and regard $\varepsilon \in \R$ as $\arg \varepsilon = \arg (1-\varepsilon) = 0$.

We denote path from p to $\varepsilon$ on real axis as $(p\varepsilon)$, then according to \eqref{HLconnect},
\begin{equation*}
\Li_{p,C}(\bw;z)=\sum_{i=0}^{r} \Li_{\varepsilon,C}(a_1\cdots a_i;z)\Li_{p,(p\varepsilon)}(a_{i+1}\cdots a_r;\varepsilon).
\end{equation*}

Since $\lim_{p \to 0} \Li_{p,(p\varepsilon)}(a_{i+1}\cdots a_r;\varepsilon) = \Li(a_{i+1}\cdots a_r;\varepsilon)$,
\begin{equation*}
\Li(\bw;z)=\sum_{i=0}^{r} \Li_{\varepsilon,C}(a_1\cdots a_i;z)\Li(a_{i+1}\cdots a_r;\varepsilon).
\end{equation*}

Consequently
\begin{align*}
|\Li(\bw;z)|&<\sum_{i=0}^{r} |\Li_{\varepsilon,C}(a_1\cdots a_i;z)||\Li(a_{i+1}\cdots a_r;\varepsilon)|\\
&<\sum_{i=0}^{r} \frac{1}{i!}\left(\frac{l(K;\varepsilon)}{\varepsilon}\right)^i<\sum_{i=0}^{\infty} \frac{1}{i!}\left(\frac{l(K;\varepsilon)}{\varepsilon}\right)^i = \exp\left(\frac{l(K;\varepsilon)}{\varepsilon}\right)
\end{align*}

In addition if $\bw \in \fH^1$,
\begin{align*}
|\Li(\bw x^n;z)| &< \sum_{j=0}^n |\Li(\reg^1(\bw x^{n-j});z)|\frac{|\log z|^j}{j!}\\
& < \exp\left(\frac{l(K;\varepsilon)}{\varepsilon}\right) \sum_{j=0}^n \frac{\ds \left(\max_{z \in K}|\log z|\right)^j}{j!} < \exp\left(\frac{l(K;\varepsilon)}{\varepsilon}\right) \exp\left(\max_{z \in K}|\log z|\right).
\end{align*}

Then let $\ds M_K=\exp\left(\frac{l(K;\varepsilon)}{\varepsilon} + \max_{z \in K}|\log z|\right)$, we get the lemma's estimation.

\end{proof}

\section{The Gauss hypergeometric function and multiple polylogarithms} \label{sec:HGF}

In this section, we consider the Gauss hypergeometric function which is the solution of the hypergeometric equation regular at the origin.

First we describe the hypergeometric functions using multiple polylogarithms as main theorem of this paper.

Next we compute the solution regular at $z=1$ and substitute these results for connection relation \eqref{connection01}'s $(1,1)$ component and get some relations of multiple zeta values.

\subsection{Representation of the hypergeometric function by multiple polylogarithms} \label{subsec:HGF:MPL}

As our first result, we construct the solution of the hypergeometric equation \eqref{HGeq} which is regular at the origin by successive iteration and write it by multiple polylogarithms.

\begin{thm}

Let parameters $\lambda_1 = \beta,\quad \lambda_2 = \alpha+1-\gamma,\quad \lambda_3 = -\alpha$.

When $|\lambda_1|,|\lambda_2|,|\lambda_3|<\frac{1}{4}$,
\begin{align}
F(\alpha,\beta,\gamma;z) &= 1 - \sum_{\substack{l,m,n \\ l,n \ge 1}} \sum_{p,q=0}^m a_{p,q}^{(l,m,n)}G_0(l+m+n,l+p,q+1;z)\lambda_1^l \lambda_2^m \lambda_3^n \label{MainThm1} \\
\frac{1}{\beta}z\frac{d}{dz}F(\alpha,\beta,\gamma;z) &= - \sum_{\substack{l,m,n \\ l,n \ge 1}} \sum_{p,q=0}^m a_{p,q}^{(l,m,n)}G(l+m+n-1,l+p,q+1;z)\lambda_1^{l-1} \lambda_2^m \lambda_3^n
\end{align}
Now,
\begin{align}
G_0(k,n,s;z) &= \sum_{\substack{k_1,\ldots,k_n \in \N\\ k_1 \ge 2 \\ k_1+\cdots+k_n=k\\\#\{i|k_i\ge 2\}=s}} \Li_{k_1,\ldots,k_n}(z) = \sum_{\substack{\bw \in \fH\\|\bw|=k-2,d(\bw)=n-1,\\s(\bw)=s-1}}\Li(x\bw y;z) \\
\intertext{i.e. Sum of all multiple polylogarithms of {\rm weight} $k$, {\rm depth} $n$,{\rm height} $s$ and $k_1 \ge 2$}
G(k,n,s;z)&= \sum_{\substack{\bw \in \fH\\|\bw|=k-1,d(\bw)=n-1,\\s(\bw)=s-1}}\Li(\bw y;z) = z\frac{d}{dz}G_0(k+1,n,s;z).
\end{align}
and coefficients
\begin{equation*}
a_{p,q}^{(l,m,n)}= \sum_{k=0}^{l-1}\binom{q}{k}\binom{l+p-q-1}{l-k-1}\binom{m+n-p-k-1}{n-1}.
\end{equation*}
Each binomial coefficients $\binom{p}{q} = 0$ if $p<0,q<0 \text{ or } p<q$.

The series of right hand side converges uniformly in the wider sense on $z \in \fU = \text{universal covering of } \bP^1-\{0,1,\infty\}$ and the representation is proper $\forall z \in \fU$.

\end{thm}

\paragraph{}

In the following, we prove the theorem.

According to Aomoto-Kita\cite{AoK}, let paramaters
\begin{equation}
\lambda_1 = \beta,\quad \lambda_2 = \alpha+1-\gamma,\quad \lambda_3 = -\alpha \label{param}
\end{equation}
and change variables as $\ds v_1 = w,\; v_2 = \frac{1}{\beta}z\frac{dw}{dz}$, \eqref{HGeq} is transformed to
\begin{equation}
\frac{d}{dz}\vect{v_1}{v_2} = (\lambda_1 \theta_1 + \lambda_2 \theta_2 + \lambda_3 \theta_3)\vect{v_1}{v_2} \label{HGeq2}
\end{equation}
\begin{align*}
\theta_1 &= \frac{1}{z}\begin{pmatrix}0 & 1 \\ 0 & 0\end{pmatrix} + \frac{1}{1-z}\begin{pmatrix}0 & 0 \\ 0 & 1\end{pmatrix} \\
\theta_2 &= \frac{1}{z}\begin{pmatrix}0 & 0 \\ 0 & 1\end{pmatrix} + \frac{1}{1-z}\begin{pmatrix}0 & 0 \\ 0 & 1\end{pmatrix} \\
\theta_3 &= \frac{1}{z}\begin{pmatrix}0 & 0 \\ 0 & 1\end{pmatrix} + \frac{1}{1-z}\begin{pmatrix}0 & 0 \\ -1 & 0\end{pmatrix}.
\end{align*}

As $\ds \vect{F(\alpha,\beta,\gamma;z)}{\frac{1}{\beta}z\frac{d}{dz}F(\alpha,\beta,\gamma;z)}$ is the solution of \eqref{HGeq2} passing $\ds \vect{1}{0}$ at $z=0$, we solve the solution by method of successive integral starting from $\vect{1}{0}$.

Namely we set vector $w_0=\vect{1}{0}$ first and determine the sequence of vectors\\
$w_n=\vect{1}{0}+\ds\int_0^z (\lambda_1 \theta_1 + \lambda_2 \theta_2 + \lambda_3 \theta_3)w_{n-1} dz$ inductively, then $\ds \lim_{n \to \infty}w_n$ is the solution of \eqref{HGeq2} passing $\vect{1}{0}$ at $z=0$.

For $\bmu = $finite sequence of $\{1,2,3\}$(even if $\bmu = \emptyset$), we define analytic functions $L_{\bmu}(z),L'_{\bmu}(z)$:
\begin{align}
\frac{d}{dz}\vect{L_{i,\bmu}(z)}{L'_{i,\bmu}(z)} &= \theta_i \vect{L_{\bmu}(z)}{L'_{\bmu}(z)} \qquad \vect{L_{i,\bmu}(0)}{L'_{i,\bmu}(0)} = \vect{0}{0} \qquad (i=1,2,3) \label{def_L} \\
\vect{L_{\emptyset}(z)}{L'_{\emptyset}(z)} &= \vect{1}{0}. \notag
\end{align}

Hence the solution can write
\begin{equation}
\vect{F(\alpha,\beta,\gamma;z)}{\frac{1}{\beta}z\frac{d}{dz}F(\alpha,\beta,\gamma;z)} = \sum_{r=1}^{\infty} \sum_{\mu_1,\ldots,\mu_r = 1,2,3} \lambda_{\mu_1}\cdots\lambda_{\mu_r}\vect{L_{\mu_1,\ldots,\mu_r}(z)}{L'_{\mu_1,\ldots,\mu_r}(z)} + \vect{1}{0}. \label{expantion1}
\end{equation}

We first express the functions $L_{\bmu}(z)$ by multiple polylogarithms. Cleary for all $\bmu$,
\begin{align*}
L_{\bmu,1}(z) &= L'_{\bmu,1}(z) = L_{\bmu,2}(z) = L'_{\bmu,2}(z) = 0 \\
L_{2,\bmu}(z) &= L_{3,\bmu}(z) = 0,
\end{align*}
we may consider the sequence only started with 1,ended by 3 for $L$ and ended with 3 for $L'$.

\begin{lem}

Define the transformation $T_0:\{\text{finite sequences of }\{1,2,3\}\} \to \fH=\Q\langle x,y\rangle$ as:
\begin{enumerate}
\item $T_0(\emptyset)=1$
\item $T_0(3,1,\bmu)=(xy-yx)T_0(\bmu)$
\item $T_0(1,\bmu)=yT_0(\bmu)$
\item $T_0(2,\bmu)=(x+y)T_0(\bmu)$
\item $T_0(3,\bmu)=xT_0(\bmu) \qquad$ if $\bk$ doesn't start with $1$
\end{enumerate}
then the functions $L_{1,\bmu,3}(z)$ and $L'_{\bmu,3}(z)$ denoted by
\begin{align}
L_{1,\bmu,3}(z) &= -\Li(xT_0(\bmu)y;z) \\
L'_{\bmu,3}(z) &= -\Li(T_0(\bmu)y;z).
\end{align}

\end{lem}

\begin{proof}

Since $\Li_{*}(z)$ and $L_{*}(z)$ are both holomorphic functions which have $0$ as $z=0$, if their derivations are identical, the original functions are so.

Also by \eqref{def_L},
\begin{align}
\frac{d}{dz}\vect{L_{1,\bmu}(z)}{L'_{1,\bmu}(z)} &= \vect{\frac{1}{z}L'_{\bmu}(z)}{\frac{1}{1-z}L'_{\bmu}(z)} \notag \\
\frac{d}{dz}\vect{L_{2,\bmu}(z)}{L'_{2,\bmu}(z)} &= \vect{0}{\frac{1}{z}L'_{\bmu}(z) + \frac{1}{1-z}L'_{\bmu}(z)} \notag \\
\frac{d}{dz}\vect{L_{3,\bmu}(z)}{L'_{3,\bmu}(z)} &= \vect{0}{\frac{1}{z}L'_{\bmu}(z) - \frac{1}{1-z}L_{\bmu}(z)}.
\end{align}

First if $L'_{\bmu,3}(z)=-\Li(T_0(\bmu)y;z)$, since
\begin{equation*}
\frac{d}{dz}L_{1,\bmu,3}(z) = \frac{1}{z}L'_{\bmu,3}(z) = -\frac{1}{z}\Li(T_0(\bmu)y;z) = -\frac{d}{dz}\Li(xT_0(\bmu)y;z)
\end{equation*}
the lemma follows also on $L$. In the following, we show $L'_{\bmu,3}(z)=-\Li(T_0(\bmu)y;z)$ by induction on length of $\bmu$.

In case of $\bmu=\emptyset$

\begin{equation*}
L'_{3}(z)=\int_0^z \frac{1}{z}L'_{\emptyset}(z)dz - \int_0^z \frac{1}{1-z}L_{\emptyset}(z)dz = - \int_0^z \frac{1}{1-z} dz = -\Li(y;z)
\end{equation*}

In Generic case, we assume $T_0(\bmu,3)=\bw y$ and $L'_{\bmu,3}(z)=-\Li(\bw y;z)$. On this time,
\begin{align*}
\frac{d}{dz}L'_{1,\bmu,3}(z) &= \frac{1}{1-z}L'_{\bmu,3}(z)=-\frac{1}{1-z}\Li(\bw y;z)=-\frac{d}{dz}\Li(y\bw y;z) \\
\frac{d}{dz}L'_{2,\bmu,3}(z) &= \frac{1}{z}L'_{\bmu,3}(z) + \frac{1}{1-z}L'_{\bmu,3}(z)=-\frac{1}{z}\Li(\bw y;z) - \frac{1}{1-z}\Li(\bw y;z) \\
&=-\frac{d}{dz}\Li(x\bw y;z) - \frac{d}{dz}\Li(y\bw y;z) = -\frac{d}{dz}\Li((x+y)\bw y;z)
\end{align*}
then,
\begin{align*}
L'_{1,\bmu,3}(z) &= -\Li(y\bw y;z) = -\Li(T_0(1,\bmu) y;z) \\
L'_{2,\bmu,3}(z) &= -\Li((x+y)\bw y;z) = -\Li(T_0(2,\bmu) y;z).
\end{align*}

Also if $\bmu$ doesn't start with $1$, then $L_{\bmu,3}(z)=0$,
\begin{equation*}
\frac{d}{dz}L'_{3,\bmu,3}(z) = \frac{1}{z}L'_{\bmu,3}(z) - \frac{1}{1-z}L_{\bmu,3}(z)=-\frac{1}{z}\Li(\bw y;z)=-\frac{d}{dz}\Li(x\bw y;z)
\end{equation*}
\begin{equation*}
L'_{3,\bmu,3}(z) = -\Li(x\bw y;z) = -\Li(T_0(3,\bmu) y;z).
\end{equation*}

Finally, since
\begin{equation*}
L_{1,\bmu,3}(z)=-\Li(x\bw y;z), \quad L'_{1,\bmu,3}(z)=-\Li(y\bw y;z)
\end{equation*}
then
\begin{align*}
\frac{d}{dz}L'_{3,1,\bmu,3}(z) &= \frac{1}{z}L'_{1,\bmu,3}(z) - \frac{1}{1-z}L_{1,\bmu,3}(z)= -\frac{1}{z}\Li(y\bw y;z) + \frac{1}{1-z}\Li(x\bw y;z) \\
&= -\frac{d}{dz}\Li(xy\bw y;z) + \frac{d}{dz}\Li(yx\bw y;z) = -\frac{d}{dz}\Li((xy-yx)\bw y;z) = -\frac{d}{dz}\Li(T_0(3,1,\bmu) y;z)
\end{align*}
Therefore
\begin{equation*}
L'_{3,1,\bmu,3}(z) = -\Li(T'_0(3,1,\bmu) y;z)
\end{equation*}

\end{proof}

By this lemma, We can rewrite \eqref{expantion1} to
\begin{align}
\vect{F(\alpha,\beta,\gamma;1)}{\frac{1}{\beta}z\frac{d}{dz}F(\alpha,\beta,\gamma;1)} &= \vect{\displaystyle \sum_{r=2}^{\infty} \sum_{\mu_2,\ldots,\mu_{r-1} = 1,2,3} \lambda_1 \lambda_{\mu_2}\cdots\lambda_{\mu_{r-1}} \lambda_3  L_{1,\mu_2,\ldots,\mu_{r-1},3}(z)}{\displaystyle \sum_{r=1}^{\infty} \sum_{\mu_1,\ldots,\mu_{r-1} = 1,2,3} \lambda_{\mu_1}\cdots\lambda_{\mu_{r-1}} \lambda_3 L'_{\mu_1,\ldots,\mu_{r-1},3}(z)} + \vect{1}{0} \label{expantion2} \notag \\
&= \vect{1}{0} -\vect{\displaystyle \sum_{r=2}^{\infty} \sum_{\mu_2,\ldots,\mu_{r-1} = 1,2,3} \lambda_1 \lambda_{\mu_2}\cdots\lambda_{\mu_{r-1}} \lambda_3 \Li(xT_0(k_2,\ldots,k_{r-1})y;z)}{\displaystyle \sum_{r=1}^{\infty} \sum_{\mu_1,\ldots,\mu_{r-1} = 1,2,3} \lambda_{\mu_1}\cdots\lambda_{\mu_{r-1}} \lambda_3 \Li(T_0(k_1,\ldots,k_{r-1})y;z)} \notag \\
&=: \sum_{l,m,n=0}^{\infty} \vect{Q(l,m,n;z)}{Q'(l,m,n;z)}\lambda_1^l \lambda_2^m \lambda_3^n.
\end{align}
i.e. $Q(0,0,0;z)=1,\; Q'(0,0,0;z)=0$ and if $(l,m,n)\neq(0,0,0)$,
\begin{align*}
Q(0,m,n;z) &= 0 \\
Q(l,m,0;z) &= 0 \\
Q'(l,m,0;z) &= 0 \\
Q(l,m,n;z) &= - \sum_{\bmu \in J(l-1,m,n-1)}\Li(xT_0(\bmu)y) \\ 
Q'(l,m,n;z) &= - \sum_{\bmu \in J(l,m,n-1)}\Li(T'_0(\bmu)y) 
\end{align*}
The index set of summation $J(l,m,n):=\{\bmu:$ {\rm finite sequence of }$\{1,2,3\}|$ {\rm $1$,$2$ and $3$ appears $l,m$ and $n$ times in $\bk$ respectively}$\}$

\paragraph{}

In this notation, if we show
\begin{equation}
Q(l,m,n;z)=-\sum_{p,q=0}^m a_{p,q}^{(l,m,n)}G_0(l+m+n,l+p,q+1;z),
\end{equation}
the theorem will been proved.

\begin{lem}

\begin{enumerate}
\item
\begin{equation}
\sum_{\bmu \in J(l,0,n)}T_0(\bmu)=\overbrace{x \cdots x}^{n\text{ \rm times}} \overbrace{y \cdots y}^{l\text{ \rm times}}
\end{equation}

Especially $Q(l,0,n;z)=-\Li_{n+1,\underbrace{\scriptstyle 1,\ldots,1}_{l-1\text{ \rm times}}}(z)$

\item

\begin{equation}
\sum_{\bmu \in J(l,m,n)}T_0(\bmu) = \sum_{p,q=0}^m \sum_{k=0}^{l}\binom{q}{k}\binom{l+p-q}{l-k}\binom{m+n-p-k}{n}
 \left( \sum_{\substack{|\bw|=l+m+n,\\d(\bw)=l+p,h(\bw)=q}}\bw \right)
\end{equation}

\end{enumerate}

\end{lem}

\begin{proof}

\begin{enumerate}

\item

Using induction on $l$. In the case of $l=0$,
\begin{equation*}
\sum_{\bmu \in J(0,0,n)}T_0(\bmu)=T_0(\underbrace{3,\ldots,3}_{n\text{ times}})=x^n.
\end{equation*}

In general cases, we assume the lemma is satisfied by $l-1$. We divide the summation by position of '1' appeared at first on the left side,
\begin{align*}
\sum_{\bmu \in J(l,0,n)}T_0(\bmu) &= \sum_{i=0}^n \sum_{\bmu \in J(l-1,0,n-i)}T_0(\underbrace{3,\ldots,3}_{i\text{ times}},1)T_0(\bmu) \\
&= yx^n y^{l-1} + \sum_{i=1}^n x^{i-1}(xy-yx)x^{n-i}y^{l-1} \\
&= yx^n y^{l-1} + \sum_{i=1}^n x^i yx^{n-i}y^{l-1} - \sum_{i=1}^n x^{i-1}yx^{n-i+1}y^{l-1}\\
&= x^n y y^{l-1} = x^n y^l.
\end{align*}

\item

We define transformation $T_0'$ by $T_0'(1)=y,T_0'(2)=x+y,T_0'(3)=x$ and $T_0'(\bmu_1,\bmu_2)=T_0'(\bmu_1)T_0'(\bmu_2)$ for $\forall \bmu_1,\bmu_2$: sequenses of $\{1,2,3\}$ and $J'(l,m,n)=\{\bmu \in J(l,m,n)|\bmu \text{ has no $(1,3)$ as subsequence}\}$.

According to (i),
\begin{equation*}
\sum_{\bmu \in J(l,0,n)}T_0(\bmu)=\overbrace{x \cdots x}^{n\text{ times}} \overbrace{y \cdots y}^{l\text{ times}} = T'_0(\overbrace{3 \cdots 3}^{n\text{ times}} \overbrace{1 \cdots 1}^{l\text{ times}}) = \sum_{\bmu \in J'(l,0,n)}T_0'(\bmu).
\end{equation*}

Therefore
\begin{equation*}
\sum_{\bmu \in J(l,m,n)}T_0(\bmu)=\sum_{\bmu \in J'(l,m,n)}T_0'(\bmu)
\end{equation*}

For $\bmu \in \text{sequence of }\{1,2,3\}$, $T_0'(\bmu)$ is $\Z$ linear combination of monomials in $\fH = \Q\langle x,y\rangle$ and if each $T_0'(\bmu)$ has component $\bw$ as summand, its multiplicity is always 1. Consequently, by fixed $\bw: |\bw|=l+m+n, d(\bw)=l+p, h(\bw)=q$, we may count number of $\{\bmu \in J'(l,m,n)| T_0'(\bmu) \text{ contains } \bw \text{ as summand}\}$ in order to compute $T_0'(\bmu)$.

The number is equal to combinatrial numbers of putting '$1$' on position of '$y$' in $\bw$, '$3$' on position of '$x$' in $\bw$ and having no juxtapoint '$13$'.
Therefore considering the situation that $k$ '1's are in position of '$y$' in '$\cdots yx \cdots$', the number can express
\begin{equation*}
\sum_{k=0}^{l}\binom{q}{k}\binom{l+p-q}{l-k}\binom{m+n-p-k}{n}.
\end{equation*}

Hence we get
\begin{equation*}
\sum_{\bmu \in J(l,m,n)}T_0(\bmu) = \sum_{p,q=0}^m \sum_{k=0}^{l}\binom{q}{k}\binom{l+p-q}{l-k}\binom{m+n-p-k}{n} \left( \sum_{\substack{|\bw|=l+m+n,\\d(\bw)=l+p,h(\bw)=q}}\bw \right).
\end{equation*}

\end{enumerate}

\end{proof}

\paragraph{Proof of theorem}

\begin{align*}
Q(l,m,n;z) &= - \sum_{\bmu \in J(l-1,m,n-1)}\Li(xT_0(\bmu)y) = - \sum_{\bmu \in J'(l-1,m,n-1)}\Li(xT'_0(\bmu)y;z) \\
&= - \sum_{p,q=0}^m \sum_{k=0}^{l-1}\binom{q}{k}\binom{l+p-q-1}{l-k-1}\binom{m+n-p-k-1}{n-1}  \sum_{\substack{|\bw|=l+m+n-2,\\d(\bw)=l+p-1,h(\bw)=q}}\Li(x\bw y ; z )
\end{align*}

Since
\begin{equation*}
\sum_{\substack{|\bw|=l+m+n-2,\\d(\bw)=l+p-1,h(\bw)=q}}\Li(x\bw y ; z ) = G_0(l+m+n,l+p,1+q;z),
\end{equation*}
\begin{align*}
Q(l,m,n;z) &= - \sum_{p,q=0}^m \sum_{k=0}^{l-1}\binom{q}{k}\binom{l+p-q-1}{l-k-1}\binom{m+n-p-k-1}{n-1} G_0(l+m+n,l+p,1+q;z) \\
&= - \sum_{p,q=0}^m a_{p,q}^{(l,m,n)} G_0(l+m+n,l+p,1+q;z).
\end{align*}

\paragraph{convergence of RHS}

First we remark that
\begin{align*}
|a_{p,q}^{(l,m,n)}| &= |\sum_{k=0}^{l-1}\binom{q}{k}\binom{l+p-q-1}{l-k-1}\binom{m+n-p-k-1}{n-1}| < \sum_{k=0}^{l-1} 2^q 2^{l+p-q-1} 2^{m+n-p-k-1} \\
&= 2^{l+m+n} \sum_{k=0}^{l-1} 2^{-k-2} < 2^{l+m+n}
\end{align*}
and $\sum_{p,q=0}^m |G_0(l+m+n,l+p,q+1;z)| < 2^{l+m+n}M_K$, because a number of indexes which have weight $k$ is $2^{k-2}<2^{k}$. 

Let $K \subset \fU$ : compact subset and assume $z \in K$. Then
\begin{align*}
&|\sum_{\substack{l,m,n \\ l,n \ge 1}} \sum_{p,q=0}^m a_{p,q}^{(l,m,n)}G_0(l+m+n,l+p,q+1;z)\lambda_1^l \lambda_2^m \lambda_3^n| \\
&\le \sum_{\substack{l,m,n \\ l,n \ge 1}} \sum_{p,q=0}^m a_{p,q}^{(l,m,n)}|G_0(l+m+n,l+p,q+1;z)||\lambda_1|^l |\lambda_2|^m |\lambda_3|^n \\
&\le \sum_{\substack{l,m,n \\ l,n \ge 1}} 2^{l+m+n} 2^{l+m+n}M_K|\lambda_1|^l |\lambda_2|^m |\lambda_3|^n \\
&\le M_K \sum_{\substack{l,m,n \\ l,n \ge 1}} |4\lambda_1|^l |4\lambda_2|^m |4\lambda_3|^n
\end{align*}

Therefore the seriese of RHS converges on $|\lambda_1|,|\lambda_2|,|\lambda_3|<\frac{1}{4}$.

\qed

\begin{cor} \label{cor:MainThm1_Cor}
\begin{align}
F(\alpha,\beta,\gamma;z)&=1-\lambda_1\lambda_3\sum_{\substack{k,n,s \in \N\\ k \ge n+s\\ n \ge s,\; s \ge 1}}G_0(k,n,s;z)(\lambda_2+\lambda_3)^{k-n-s}(\lambda_1+\lambda_2)^{n-s}\lambda_2^{s-1}(\lambda_1+\lambda_2+\lambda_3)^{s-1} \label{MainThm1_Cor} \\
\frac{1}{\beta}z\frac{d}{dz}F(\alpha,\beta,\gamma;z)&=- \lambda_3  \sum_{k,n,s} G(k-1,n,s;z) (\lambda_2+\lambda_3)^{k-n-s} (\lambda_1+\lambda_2)^{n-s} \lambda_2^{s-1}(\lambda_1+\lambda_2+\lambda_3)^{s-1} \label{MainThm1_Cor2}
\end{align}
and convergent condition can be expand to $|\lambda_1+\lambda_2|,|\lambda_2+\lambda_3|,|\lambda_2|,|\lambda_1+\lambda_2+\lambda_3|<1$.

\end{cor}

\begin{proof}

We remark $\binom{(m-p)+n-k}{n}=\sum_{t=0}^n \binom{(m-p)+n-q}{n-t}\binom{q-k}{t} \quad k \le \forall q \le (m-p)+n$.

\begin{align*}
&\sum_{\substack{l,m,n\\ l,n \ge 1}}\sum_{p,q=0}^m a_{p,q}^{(l,m,n)}G_0(l+m+n,l+p,q+1;z)\lambda_1^l\lambda_2^m\lambda_3^n \\
&= \lambda_1\lambda_3\sum_{l,m,n \in \Z_{\ge 0}}\sum_{p,q=0}^m \sum_{k=0}^l \binom{q}{k}\binom{l+p-q}{l-k}\binom{m+n-p-k}{n}G_0(l+m+n+2,l+p+1,q+1;z)\lambda_1^l\lambda_2^m\lambda_3^n \\
&= \lambda_1\lambda_3\sum_{l,m,n \in \Z_{\ge 0}}\sum_{p,q=0}^m \sum_{k=0}^l \sum_{t=0}^n \binom{(m-p)+n-q}{n-t}\binom{l+p-q}{l-k}\binom{q}{k}\binom{q-k}{t}\\
&\phantom{ \lambda_1\lambda_3\sum_{l,m,n \in \Z_{\ge 0}}\sum_{p,q=0}^m \sum_{k=0}^l \sum_{t=0}^n}\times G_0(l+m+n+2,l+p+1,q+1;z)\lambda_2^{m-p-q+t}\lambda_3^{n-t}\lambda_1^{l-k}\lambda_2^{p-q+k} \lambda_2^q \lambda_1^k\lambda_2^{q-k-t}\lambda_3^t \\
\intertext{Now let $k'=m+n-p-q, n'=l+p-q, s'=q$,}
&= \lambda_1\lambda_3\sum_{k',n',s' \in \Z_{\ge 0}}G_0(k'+n'+2s'+2,n'+s'+1,s'+1;z)(\lambda_2+\lambda_3)^{k'}(\lambda_1+\lambda_2)^{n'} \lambda_2^{s'} (\lambda_1+\lambda_2+\lambda_3)^{s'} \\
&= \lambda_1\lambda_3\sum_{\substack{k,n,s \in \N\\ k \ge n+s\\ n \ge s,\; s \ge 1}}G_0(k,n,s;z)(\lambda_2+\lambda_3)^{k-n-s}(\lambda_1+\lambda_2)^{n-s} \lambda_2^{s-1} (\lambda_1+\lambda_2+\lambda_3)^{s-1}
\end{align*}

The claim about $\frac{1}{\beta}z\frac{d}{dz}F$ and convergent condition are clear.

\end{proof}

\paragraph{Examples of theorem}

In cases of some lower $m$, (and $\forall \; l,n \ge 1$)
\begin{align}
\text{coefficient of $\lambda_1^l \lambda_3^n$} &= - \Li_{n+1,\underbrace{\scriptstyle 1,\ldots,1}_{l-1\text{ times}}}(z) \quad (=G_0(l+n,l,1;z))\\
\text{coefficient of $\lambda_1^l \lambda_2 \lambda_3^n$} &= - n G_0(l+n+1,l,1;z) - G_0(l+n+1,l,2;z) \\
&\qquad - l G_0(l+n+1,l+1,1;z) - G_0(l+n+1,l+1,2;z) \notag \\
\text{coefficient of $\lambda_1^l \lambda_2^2 \lambda_3^n$} &= -\frac{n(n+1)}{2} G_0(l+n+2,l,1;z) - n G_0(l+n+2,l,2;z) \notag \\
&- G_0(l+n+2,l,3;z) - ln G_0(l+n+2,l+1,1;z) \notag \\
&- (l+n-1) G_0(l+n+2,l+1,2;z) - 2 G_0(l+n+2,l+1,3;z) \notag \\
&- \frac{l(l+1)}{2} G_0(l+n+2,l+2,1;z) - l G_0(l+n+2,l+2,2;z) \notag \\
&- G_0(l+n+2,l+2,3;z).
\end{align}

Consequently, lower degree terms of $F(\alpha,\beta,\gamma;z)$ in $\lambda_1,\lambda_2$ and $\lambda_3$ is 
\begin{align}
F(\alpha,\beta,\gamma;z) &= 1 - \Li_{2}(z)\lambda_1\lambda_3 - \Li_{2,1}(z)\lambda_1^2\lambda_3 - \Li_{3}(z)\lambda_1\lambda_3^2 - (\Li_{3}(z)+\Li_{2,1}(z))\lambda_1\lambda_2\lambda_3 \notag \\
& - \Li_{2,1,1}(z)\lambda_1^3\lambda_3 - \Li_{3,1}(z)\lambda_1^2\lambda_3^2 - \Li_{4}(z)\lambda_1\lambda_3^3 - (\Li_{3,1}(z)+\Li_{2,2}(z)+2\Li_{2,1,1}(z))\lambda_1^2\lambda_2\lambda_3 \notag \\
& - (\Li_{4}(z)+\Li_{3,1}(z)+\Li_{2,2}(z)+\Li_{2,1,1}(z))\lambda_1\lambda_2^2\lambda_3 - (2\Li_{4}(z)+\Li_{3,1}(z)+\Li_{2,2}(z))\lambda_1\lambda_2\lambda_3^2 \notag \\
& - \cdots
\end{align}

\subsection{Regular solutions in neighborhood of $1$} \label{subsec:HGF:nbh1}

Next, we construct the regular solution of hypergeometric equation in neighborhood of $1$.

We consider ${}^t \Phi_1^{-1}(1-z)$. Let $t = 1-z$.

Since $\frac{d}{dt}\Phi_1^{-1}(t) = -\Phi_1^{-1} \frac{d \Phi_1}{dt}\Phi_1^{-1}$, ${}^t \Phi_1^{-1}(t)$ satisfies the equation
\begin{align}
\frac{d}{dt} {}^t \Phi_1^{-1}(t) &= \left( \frac{1}{t}\begin{pmatrix}0&\alpha \\ 0 & \alpha+\beta+1-\gamma\end{pmatrix} + \frac{1}{1-t}\begin{pmatrix}0 & 0 \\ \beta & 1-\gamma\end{pmatrix} \right) {}^t \Phi_1^{-1}(t) \notag \\
 &= (\lambda_1 \theta'_1 + \lambda_2 \theta'_2 + \lambda_3 \theta'_3) {}^t \Phi_1^{-1}(t)\label{Phi_inv_eq}
\end{align}
\begin{align*}
\theta'_1 &= \frac{1}{t}\begin{pmatrix}0 & 0 \\ 0 & 1\end{pmatrix} + \frac{1}{1-t}\begin{pmatrix}0 & 0 \\ 1 & 0\end{pmatrix} \\
\theta'_2 &= \frac{1}{t}\begin{pmatrix}0 & 0 \\ 0 & 1\end{pmatrix} + \frac{1}{1-t}\begin{pmatrix}0 & 0 \\ 0 & 1\end{pmatrix} \\
\theta'_3 &= \frac{1}{t}\begin{pmatrix}0 & -1 \\ 0 & 0\end{pmatrix} + \frac{1}{1-t}\begin{pmatrix}0 & 0 \\ 0 & 1\end{pmatrix}
\end{align*}

We denote
\begin{equation}
\Phi^{-1}_1(t)=\begin{pmatrix}\psi_{11}(t) & \psi_{12}(t) \\ \psi_{21}(t) & \psi_{22}(t)\end{pmatrix}.
\end{equation}

$\psi_{11}(t)$ and $\psi_{12}(t)$ are regular at $t=0$ and $\psi_{11}(0)=1$,$\psi_{12}(0)=0$ for direct culculation.

Therefore the solution of the equation which is regular at $t=0$: $\ds \vect{\psi_{11}(t)}{\psi_{12}(t)}$ can be solved in the same way as section \ref{sec:HGF}.

We express
\begin{equation}
\vect{\psi_{11}(t)}{\psi_{12}(t)} = \sum_{r=1}^{\infty} \sum_{\mu_1,\ldots,\mu_r = 1,2,3} \lambda_{\mu_1}\cdots\lambda_{\mu_r}\vect{L^{(1)}_{\mu_1,\ldots,\mu_r}(t)}{L^{'(1)}_{\mu_1,\ldots,\mu_r}(t)} + \vect{1}{0},
\end{equation}
then functions $L^{(1)}_{\mu_1,\ldots,\mu_r}(t),L^{'(1)}_{\mu_1,\ldots,\mu_r}(t)$ are determined by the following lemma.

\begin{lem}

Define the transformation $T_1:\{\text{finite sequences of}\{1,2,3\}\} \to \fH$ as:
\begin{enumerate}
\item $T_1(\emptyset)=1$
\item $T_1(1,3,\bmu)=(xy-yx)T_1(\bmu)$
\item $T_1(1,\bmu)=xT_1(\bmu) \qquad$ if $\bmu$ doesn't start with 3
\item $T_1(2,\bmu)=(x+y)T_1(\bmu)$
\item $T_1(3,\bmu)=yT_1(\bmu)$.
\end{enumerate}

Then the functions $L^{(1)}_{3,\bmu,1}(t),L^{'(1)}_{\bmu,1}(t)$ denoted by
\begin{align}
L^{(1)}_{3,\bmu,1}(t) &= -\Li(xT_1(\bmu)y;t) \\
L^{'(1)}_{\bmu,1}(t) &= \Li(T_1(\bmu)y;t).
\end{align}

\end{lem}

This lemma shows $T_1=\tau \circ T_0$, then $\psi_{11}(t)$ is "dual" of $\varphi^{(0)}_0(z)$.

Cleary $\ds \sum_{\bw \in I(k,n,s)}x\bw y = \sum_{\bw \in I(k,k-n,s)}x\tau(\bw) y$ and $\ds \sum_{\bw \in I(k,n,s)}\bw y = \sum_{\bw \in I(k,k-n,s)}\tau(\bw) y$, we get the following proposition.

\begin{prop}
\begin{align}
\psi_{11}(1-z) &= 1 - \sum_{\substack{l,m,n \\ l,n \ge 1}} \sum_{p,q=0}^m a_{p,q}^{(l,m,n)}G_0(l+m+n,m+n-p,q+1;1-z)\lambda_1^l \lambda_2^m \lambda_3^n \\
&= 1 - \lambda_1 \lambda_3  \sum_{k,n,s} G_0(k,k-n,s;1-z) (\lambda_2+\lambda_3)^{k-n-s} (\lambda_1+\lambda_2)^{n-s} \lambda_2^{s-1}(\lambda_1+\lambda_2+\lambda_3)^{s-1} \\
\psi_{12}(1-z) &= \sum_{\substack{l,m,n \\ l,n \ge 1}} \sum_{p,q=0}^m a_{p,q}^{(l,m,n)}G(l+m+n-1,m+n-p,q+1;1-z)\lambda_1^l \lambda_2^m \lambda_3^{n-1} \\
&= \lambda_1  \sum_{k,n,s} G(k-1,k-n,s;1-z) (\lambda_2+\lambda_3)^{k-n-s} (\lambda_1+\lambda_2)^{n-s} \lambda_2^{s-1}(\lambda_1+\lambda_2+\lambda_3)^{s-1}
\end{align}
\end{prop}

\subsection{The connection relation between the regular solution of $0$ and $1$}\label{subsec:HGF:connection}

We consider the $(1,1)$ component of connection relation(\ref{connection01}):
\begin{equation}
\psi_{11}(1-z)\varphi^{(0)}_0(z) + \psi_{12}(1-z) \frac{1}{\beta}z\frac{d}{dz}\varphi^{(0)}_0(z) = \frac{\ds \Gamma(\gamma)\Gamma(\gamma-\alpha-\beta)}{\ds \Gamma(\gamma-\alpha)\Gamma(\gamma-\beta)}. \label{HGFconnection_11}
\end{equation}

\begin{thm}

\begin{align}
&1-\lambda_1\lambda_3\sum_{\substack{k,n,s \in \N\\ k \ge n+s\\ n \ge s,\; s \ge 1}}(G_0(k,n,s;z)+G_0(k,k-n,s;1-z))(\lambda_2+\lambda_3)^{k-n-s}(\lambda_1+\lambda_2)^{n-s}\lambda_2^{s-1}(\lambda_1+\lambda_2+\lambda_3)^{s-1} \notag \\
&+\lambda^2_1\lambda^2_3 \sum_{k,n,s}\sum_{k'=2}^{k-2} \sum_{n'=1}^{n-1} \sum_{s'=1}^{s-1}G_0(k',n',s';z)G_0(k-k',(k-k')-(n-n'),s-s';1-z) \notag \\
&\phantom{+\lambda^2_1\lambda^2_3 \sum_{k,n,s}\sum_{k'=2}^{k-2} \sum_{n'=1}^{n-1} \sum_{s'=1}^{s-1}}   \times (\lambda_2+\lambda_3)^{k-n-s}(\lambda_1+\lambda_2)^{n-s}\lambda_2^{s-2}(\lambda_1+\lambda_2+\lambda_3)^{s-2}\notag \\
&- \lambda^2_3  \sum_{k,n,s}\sum_{k'=2}^{k-2} \sum_{n'=1}^{n-1} \sum_{s'=1}^{s-1} G(k'-1,n',s';z)G(k-k'-1,(k-k')-(n-n'),s-s';1-z) \notag\\
&\phantom{- \lambda^2_3  \sum_{k,n,s}\sum_{k'=2}^{k-2} \sum_{n'=1}^{n-1} \sum_{s'=1}^{s-1}} \times (\lambda_2+\lambda_3)^{k-n-s} (\lambda_1+\lambda_2)^{n-s} \lambda_2^{s-2}(\lambda_1+\lambda_2+\lambda_3)^{s-2} \notag\\
&= \frac{\ds \Gamma(1-(\lambda_2+\lambda_3))\Gamma(1-(\lambda_1+\lambda_2))}{\ds \Gamma(1-\lambda_2)\Gamma(1-(\lambda_1+\lambda_2+\lambda_3))} \label{MainThm2}
\end{align}

Especially let $\lambda_2=0$, we get
\begin{align}
1&-\sum_{l,n \ge 1}(G_0(l+n,l,1;z)+G_0(l+n,n,1;1-z))\lambda_1^l \lambda_3^n \notag \\
&+ \sum_{l,n \ge 2}\sum_{l'=0}^{l-1}\sum_{n'=0}^{n-1}G_0(l'+n',l',1;z)G_0(l-l'+n-n',n-n',1;1-z)\lambda_1^l \lambda_3^n \notag \\
&- \sum_{l,n \ge 1}\sum_{l'=0}^{l-1}\sum_{n'=0}^{n-1}G(n-n'+l',l'+1,1;z)G(l-l'+n',n'+1,1;1-z)\lambda_1^l \lambda_3^n \notag \\
&= \frac{\ds \Gamma(1-\lambda_3)\Gamma(1-\lambda_1)}{\ds \Gamma(1-(\lambda_1+\lambda_3))}
\end{align}

\end{thm}

\paragraph{}

This formula is generalization of Euler's inversion formula. Indeed comparing coefficients of $\lambda_1\lambda_3^n$, we get
\begin{equation}
\Li_{n+1}(z)+\Li_{2,\underbrace{\scriptstyle 1,\ldots,1}_{n-1\text{ times}}}(1-z)+\sum_{j=1}^n\Li_{n-j+1}(z)\Li_{\underbrace{\scriptstyle 1,\ldots,1}_{j\text{ times}}}(1-z)=\zeta(n+1).
\end{equation}
The special case of $n=1$ as this formula is classical Euler's inversion formula for di-logarighm.

\paragraph{}

We can get many relations of multiple zeta values as limit of \eqref{MainThm2} as $z \to 0,1$.

Now,
\begin{align*}
G_0(k,n,s;z) &\to G_0(k,n,s;1) \qquad z \to 1\\
G_0(k,n,s;1-z) &\to 0 \qquad z \to 1 \text{ as polynomial order}\\
G(k,n,s;z) &\to \infty \qquad z \to 1 \text{ as logarithm order}\\
G(k,n,s;1-z) &\to 0 \qquad z \to 1 \text{ as polynomial order}
\end{align*}
then the equation of theorem \eqref{MainThm2} tend to
\begin{align}
F(\alpha,\beta,\gamma;z)&= 1-\lambda_1\lambda_3\sum_{\substack{k,n,s \in \N\\ k \ge n+s\\ n \ge s,\; s \ge 1}}G_0(k,n,s;1)(\lambda_2+\lambda_3)^{k-n-s}(\lambda_1+\lambda_2)^{n-s}\lambda_2^{s-1}(\lambda_1+\lambda_2+\lambda_3)^{s-1} \notag \\
&= 1 - \sum_{\substack{l,m,n \\ l,n \ge 1}} \sum_{p,q=0}^m a_{p,q}^{(l,m,n)}G_0(l+m+n,l+p,q+1;1)\lambda_1^l \lambda_2^m \lambda_3^n \notag \\
&= \frac{\ds \Gamma((\lambda_2+\lambda_3))\Gamma(\lambda_1+\lambda_2)}{\ds \Gamma(\lambda_2)\Gamma(\lambda_1+\lambda_2+\lambda_3)} \label{OZrel}
\end{align}
as $z \to 1$.

In the same way, \eqref{MainThm2} converge to
\begin{align}
&1-\lambda_1\lambda_3\sum_{\substack{k,n,s \in \N\\ k \ge n+s\\ n \ge s,\; s \ge 1}}G_0(k,k-n,s;1)(\lambda_2+\lambda_3)^{k-n-s}(\lambda_1+\lambda_2)^{n-s}\lambda_2^{s-1}(\lambda_1+\lambda_2+\lambda_3)^{s-1} \notag \\
&= 1 - \sum_{\substack{l,m,n \\ l,n \ge 1}} \sum_{p,q=0}^m a_{p,q}^{(l,m,n)}G_0(l+m+n,m+n-p,q+1;1)\lambda_1^l \lambda_2^m \lambda_3^n \notag \\
&= \frac{\ds \Gamma(1-(\lambda_2+\lambda_3))\Gamma(1-(\lambda_1+\lambda_2))}{\ds \Gamma(1-\lambda_2)\Gamma(1-(\lambda_1+\lambda_2+\lambda_3))}
\end{align}
as $z \to 0$. Then comparing both limits, we get the duality formula of multiple zeta value partially.

\begin{cor}[duality]
\begin{equation}
G_0(k,n,s;1)=G_0(k,k-n,s;1)
\end{equation}
\end{cor}

The duality formula is also understood symmetry between $\lambda_1$ and $\lambda_3$ in \label{OZrel}.

\paragraph{}

We can expand the RHS of \eqref{OZrel} as seriese of $\lambda_1,\lambda_2,\lambda_3$.

Next expanding gamma functions on RHS by
\begin{equation}
\Gamma(1-z) = \exp(cz - \sum_{n=2}^{\infty}\frac{\zeta(n)}{n}z^n) \qquad \qquad (\ds c=\lim_{n \to \infty} (\sum_{k=1}^n \frac{1}{k}-\log n)\quad \text{: Euler constant})
\end{equation}
and
\begin{equation*}
\exp(\sum_{k=1}^{\infty} x_k z^k) = \sum_{k=0}^{\infty} S_k(\bx) z^k \qquad \bx=(x_1,x_2,\ldots)
\end{equation*}
\begin{equation*}
S_k(\bx)=\sum_{k_1+2k_2+3k_3+\cdots = k} \frac{x_1^{k_1}}{k_1!}\frac{x_2^{k_2}}{k_2!}\frac{x_3^{k_3}}{k_3!}\cdots \qquad \text{(Schur polynomial)},
\end{equation*}
we get the next theorem.

\begin{thm}

\begin{multline}
\sum_{p,q=0}^m \left(\sum_{k=0}^{l-1}\binom{q}{k}\binom{l+p-q-1}{l-k-1}\binom{m+n-p-k-1}{n-1}\right)G_0(l+m+n,l+p,q+1;1)\\
 = \sum_{\substack{l_2+l_4=l \\ n_1+n_4=n \\ m_1+m_2+m_3+m_4=m}}\frac{(m_1+n_1)!}{m_1!n_1!}\frac{(l_2+m_2)!}{l_2!m_2!}\frac{(l_4+m_4+n_4)!}{l_4!m_4!n_4!}\\
 \times S_{m_1+n_1}({\text{\boldmath $\zeta$}})S_{l_2+m_2}({\text{\boldmath $\zeta$}})S_{m_3}(-{\text{\boldmath $\zeta$}})S_{l_4+m_4+n_4}(-{\text{\boldmath $\zeta$}}) \label{MainThm3}
\end{multline}
for $\displaystyle {\text{\boldmath $\zeta$}} = (0,\frac{\zeta(2)}{2},\frac{\zeta(3)}{3},\frac{\zeta(4)}{4},\ldots)$.

\end{thm}

\paragraph{}

In addition, regarding both side of (\ref{MainThm3}) as a polynomial in $\lambda_1$ and comparing both constant terms, we get the sum formula of multiple zeta values shown by Granville\cite{G} and Zagier\cite{Z2}.

\begin{cor}[Sum formula]
\begin{equation}
\sum_{s}G_0(w,d,s)=\zeta(w) \label{Sum_formula}
\end{equation}

\end{cor}

\begin{proof}

We divide both side by $\lambda_1\lambda_3$ and let $\lambda_1 \to 0$. We denote $\displaystyle \sum_{s}G_0(w,d,s) = G_0(w,d,*)$.

\paragraph{LHS}

\begin{align*}
\sum_{\substack{l,m,n\\ l,n\ge 1}} & \sum_{p,q=0}^m \left(\sum_{k=0}^{l-1}\binom{q}{k}\binom{l+p-q-1}{l-k-1}\binom{m+n-p-k-1}{n-1}\right)G_0(l+m+n,l+p,q+1)\lambda_1^{l-1} \lambda_2^m \lambda_3^{n-1} \\
& \overset{\lambda_1 \to 0}{\longrightarrow} \sum_{\substack{m,n\\ n\ge 1}} \sum_{p,q=0}^m \left(\binom{q}{0}\binom{p-q}{0}\binom{m+n-p-1}{n-1}\right)G_0(m+n+1,p+1,q+1)\lambda_2^m \lambda_3^{n-1} \\
&= \sum_{m,n} \sum_{p,q=0}^m \binom{m+n-p}{n}G_0(m+n+2,p+1,q+1)\lambda_2^m \lambda_3^n \\
&= \sum_{w=0}^{\infty} \sum_{n=0}^{w} \sum_{p,q=0}^{w-n} \binom{w-p}{n}G_0(w+2,p+1,q+1)\lambda_2^{w-n} \lambda_3^n \\
&= \sum_{w=0}^{\infty} \sum_{n=0}^{w} \sum_{p=0}^{w-n} \binom{w-p}{n}G_0(w+2,p+1,*)\lambda_2^{w-n} \lambda_3^n \\
&= \sum_{w=0}^{\infty} \sum_{p=0}^{w} G_0(w+2,p+1,*) \sum_{n=0}^{w-p} \binom{w-p}{n}\lambda_2^p \lambda_2^{w-p-n} \lambda_3^n \\
&= \sum_{w=0}^{\infty} \sum_{p=0}^{w} G_0(w+2,p+1,*) \lambda^p(\lambda_2+\lambda_3)^{w-p} \\
&= \sum_{w=2}^{\infty} \sum_{d=1}^{w-1} G_0(w,d,*) \lambda_2^{d-1} (\lambda_2+\lambda_3)^{w-d-1}
\end{align*}

\paragraph{RHS}

We partial differntiate RHS by $\lambda_1$ , let $\lambda_1 \to 0$ and divide by $\lambda_3$:
\begin{align*}
1-\exp&\left(\sum_{n \ge 2} \frac{\zeta(n)}{n}\left((\lambda_2+\lambda_3)^n+(\lambda_1+\lambda_2)^n-\lambda_2^n-(\lambda_1+\lambda_2+\lambda_3)^n\right)\right) \\
& \overset{\partial/\partial \lambda_1}{\longrightarrow} -\sum_{n \ge 2} \frac{\zeta(n)}{n}\left(n(\lambda_1+\lambda_2)^{n-1}-n(\lambda_1+\lambda_2+\lambda_3)^{n-1}\right) \\
& \overset{\lambda_1 \to 0}{\longrightarrow} -\sum_{n \ge 2} \zeta(n)\left(\lambda_2^{n-1}-(\lambda_2+\lambda_3)^{n-1}\right) \\
& = \sum_{n \ge 2} \zeta(n) \lambda_3 \left(\lambda_2^{n-2} + \lambda_2^{n-3}(\lambda_2+\lambda_3) + \cdots + \lambda_2(\lambda_2+\lambda_3)^{n-3} + (\lambda_2+\lambda_3)^{n-2} \right) \\
& \overset{1/\lambda_3}{\longrightarrow}\sum_{n \ge 2} \zeta(n) \sum_{k=0}^{n-2}\lambda_2^k(\lambda_2+\lambda_3)^{n-k-2} \\
& = \sum_{n \ge 2} \sum_{k=1}^{n-1} \zeta(n) \lambda_2^{k-1}(\lambda_2+\lambda_3)^{n-k-1}
\end{align*}

\paragraph{}

Therefore comparing both side's coefficients with respect to $\lambda_2$ and $(\lambda_2+\lambda_3)$, we get the equation \eqref{Sum_formula}.

\end{proof}

\paragraph{}

On the above results, Ohno-Zagier\cite{OZ} showed previously that generation function of $G_0(k,n,s;t)$:
\begin{equation}
\sum_{k,n,s}G_0(k,n,s;t)x^{k-n-s}y^{n-s}z^{s-1} \label{OZseries}
\end{equation}
satisfies the hypergeometric equation and the function is expressed as
\begin{align}
\sum_{k,n,s}G_0(k,n,s;t)x^{k-n-s}y^{n-s}z^{s-1} &= \frac{1}{xy-z}(1-F(\alpha-x,\beta-x,1-x;t)) \label{OZrelation} \\
\intertext{and}
\sum_{k,n,s}G_0(k,n,s;1)x^{k-n-s}y^{n-s}z^{s-1} &= \frac{1}{xy-z}(1-\exp(\sum_{n=2}^{\infty}\frac{\zeta(n)}{n}(x^n+y^n-\alpha^n-\beta^n)))
\end{align}
(Now,$\alpha,\beta$ is complex numbers satisfying $\alpha+\beta=x+y$, $\alpha\beta=z$).

And they showed some relations of multiple zeta values using Gauss formula (\ref{GAUSSformula}). In particular, specializing the paramater $x,y$ and $z$, they showed the sum formula, the Le-Murakami formula and some other formulas of multiple zeta values.

\eqref{OZrel} means that our results are equal to relations of Ohno-Zagier by change of variables: $x=\lambda_2+\lambda_3,y=\lambda_1+\lambda_2,z=\lambda_2(\lambda_1+\lambda_2+\lambda_3)$. Consequently the results is not new essentially, but give natural interpretation of relation of Ohno-Zagier as solving the hypergeometric equation directly and can be developed as replacing the hypergeometric equation by many variable KZ equation.

\section{On singular solutions of hypergeometric equation} \label{sec:OtherSolution}

In this section, we consider singular solutions of the hypergeometric equation in neighborhood of $0,1,\infty$. We assume $\alpha, \beta, \gamma, \gamma-\alpha-\beta \not\in \Z$ throught this section.

\subsection{Solutions singular at the origin} \label{subsec:OtherSolution:singular}

First we consider the $(1,2)$ component of \eqref{connection01}.

We now start with $\vect{0}{1}$ instead of $\vect{1}{0}$ on section \ref{sec:HGF}, then
\begin{equation}
\vect{\tilde{\varphi}(z)}{\tilde{\varphi}'(z)} = \vect{0}{1} +\vect{\displaystyle \sum_{r=0}^{\infty} \sum_{\mu_1,\ldots,\mu_r = 1,2,3} \lambda_1 \lambda_{\mu_1}\cdots\lambda_{\mu_r} \Li(xT_0(\mu_1,\ldots,\mu_r);z)}{\displaystyle \sum_{r=1}^{\infty} \sum_{\mu_1,\ldots,\mu_r = 1,2,3} \lambda_{\mu_1}\cdots\lambda_{\mu_r} \Li(T_0(\mu_1,\ldots,\mu_r);z)} \label{singular_solution_of_HGeq}
\end{equation}
is a solution of equation (\ref{HGeq2}).

We compute asymptotic condition of $\vect{\tilde{\varphi}(z)}{\tilde{\varphi}'(z)}$ as $z \to 0$.

Now $f(z) \to g(z) \quad (z \to 0)$ means $f(z)g^{-1}(z) \to 1 \quad (z \to 0)$.

By definition of $\Li$ and $\Li(\bw;z)\log^n z \to 0 \quad (z \to 0)$ for $\bw \in \fH^{1}-\Q, n \in \Z_{\ge 0}$,
\begin{equation}
\Li(\bw x^n;z) \overset{z \to 0}{\longrightarrow} \begin{cases}\frac{\log^n z}{n!},&(\bw \in \Q)\\ 0,&(\bw \not\in \Q)\end{cases}
\end{equation}
and monomial $x^n$ is generated by $xT_0(\mu_1,\ldots,\mu_r)$ only in case of $r=n-1$ and $\mu_1,\ldots,\mu_{n-1}= 2$ or $3$.

Hence asymptotic condition of $\tilde{\varphi}(z)$ is
\begin{align}
\sum_{r=0}^{\infty} \lambda_1 (\lambda_2+\lambda_3)^r \frac{\log^{r+1} z}{(r+1)!} &= \frac{\lambda_1}{\lambda_2+\lambda_3} \sum_{r=0}^{\infty} \frac{((\lambda_2+\lambda_3)\log z)^{r+1}}{(r+1)!} \notag \\
&= \frac{\lambda_1}{\lambda_2+\lambda_3} \left(\sum_{r=0}^{\infty} \frac{((\lambda_2+\lambda_3)\log z)^r}{r!} - 1 \right) \notag \\
&= \frac{\lambda_1}{\lambda_2+\lambda_3} \left(z^{\lambda_2+\lambda_3} - 1 \right).
\end{align}

In similar way, asymptotic condition of $\tilde{\varphi}'(z)$ as $z \to 0$ is $z^{\lambda_2+\lambda_3}$.

Since asymptotic condition of $\vect{\varphi^{(0)}_1(z)}{\frac{1}{\beta}z\frac{d}{dz}\varphi^{(0)}_1(z)}$ as $z \to 0$ is
\begin{equation}
\vect{\varphi^{(0)}_1(z)}{\frac{1}{\beta}z\frac{d}{dz}\varphi^{(0)}_1(z)} = \vect{z^{1-\gamma}F(\alpha-\gamma+1,\beta-\gamma+1,2-\gamma;z)}{z\frac{d}{dz}\left(z^{1-\gamma}F(\alpha-\gamma+1,\beta-\gamma+1,2-\gamma;z)\right)} \overset{z \to 0}{\longrightarrow}\vect{z^{\lambda_2+\lambda_3}}{\frac{\lambda_2+\lambda_3}{\lambda_1}z^{\lambda_2+\lambda_3}},
\end{equation}
the solution of the hypergeometric equation \eqref{HGeq2}: $\vect{\varphi^{(0)}_1(z)}{\frac{1}{\beta}z\frac{d}{dz}\varphi^{(0)}_1(z)}$ can be written as
\begin{align}
\vect{\varphi^{(0)}_1(z)}{\frac{1}{\beta}z\frac{d}{dz}\varphi^{(0)}_1(z)} = \frac{\lambda_2+\lambda_3}{\lambda_1} \vect{\tilde{\varphi}(z)}{\tilde{\varphi}'(z)} + \vect{\varphi^{(0)}_0(z)}{\frac{1}{\beta}z\frac{d}{dz}\varphi^{(0)}_0(z)}. \label{solution_singular_0}
\end{align}

Then we get the following proposition.

\begin{prop}
\begin{align}
&z^{1-\gamma}F(\alpha-\gamma+1,\beta-\gamma+1,2-\gamma;z) \notag \\
& = \sum_{l,m,n=0}^{\infty} \lambda_1^l \lambda_2^m \lambda_3^n \bigg(\sum_{\bmu \in J'(l-1,m-1,n)}\Li(xT'_0(\bmu)y;z) + \sum_{j=0}^n\sum_{\bmu \in J'(l,m-2,n-j)}\Li(xT'_0(\bmu)(x+y)x^j;z) \notag\\
& \phantom{+ \sum_{j=0}^n\sum_{\bmu \in J'(l,m-2,n-j)}\Li(xT'_0(\bmu)(x+y)x^j;z)} + \sum_{j=0}^{n-1}\sum_{\bmu \in J'(l,m-1,n-j-1)}\Li(xT'_0(\bmu)(x+y)x^j;z) \bigg) + 1
\end{align}

And as $z \to 1$,
\begin{equation}
\begin{split}
&\sum_{l,m,n=0}^{\infty} \lambda_1^l \lambda_2^m \lambda_3^n \bigg(\sum_{\bmu \in J'(l-1,m-1,n)}\Li(xT'_0(\bmu)y;1) + \sum_{j=0}^n\sum_{\bmu \in J'(l,m-2,n-j)}\Li(xT'_0(\bmu)(x+y)x^j;1) \\
& \phantom{+ \sum_{j=0}^n\sum_{\bmu \in J'(l,m-2,n-j)}} + \sum_{j=0}^{n-1}\sum_{\bmu \in J'(l,m-1,n-j-1)}\Li(xT'_0(\bmu)(x+y)x^j;1) \bigg) = \frac{\ds \Gamma(1+(\lambda_2+\lambda_3))\Gamma(1-(\lambda_1+\lambda_2))}{\ds \Gamma(1+\lambda_3)\Gamma(1-\lambda_1)} -1 \label{connection_12}
\end{split}
\end{equation}

\end{prop}

\begin{proof}

The first claim can be proved by computing \eqref{solution_singular_0} by \eqref{MainThm1} and \eqref{singular_solution_of_HGeq}.

The second claim is limit as $z \to 1$ of connection relation \eqref{connection01}'s $(1,2)$ component:
\begin{equation}
\psi_{11}(1-z)\varphi^{(0)}_1(z) + \psi_{12}(1-z) \frac{1}{\beta}z\frac{d}{dz}\varphi^{(0)}_1(z) = \frac{\ds \Gamma(2-\gamma)\Gamma(\gamma-\alpha-\beta)}{\ds \Gamma(1-\alpha)\Gamma(1-\beta)}.
\end{equation}

\end{proof}

We can compute expression \eqref{connection_12} by using following lemma.

\begin{lem} \label{lem:culculus_hH}
Assuming $\bw = x^{k_1-1}yx^{k_2-1}y\cdots x^{k_r-1}$ and $k_1 \ge 2$,
\begin{equation}
\Li(\bw yx^n;1)= (-1)^n \sum_{\varepsilon_1+\cdots+\varepsilon_r=n} \binom{k_1+\varepsilon_1-1}{\varepsilon_1}\cdots\binom{k_r+\varepsilon_r-1}{\varepsilon_r}\zeta(k_1+\varepsilon_1,\ldots,k_r+\varepsilon_r).
\end{equation}

\end{lem}

\begin{proof}

If $\bw \in x\fH$,
\begin{equation*}
\Li(\bw y x^n;1) = \Li(\reg^{1}(\bw y x^n);1) = (-1)^n \Li((\bw \sh x^n)y;1) = (-1)^n \frac{1}{n!}\Li((\bw \sh x^{\sh n})y;1).
\end{equation*}

Then we prove
\begin{equation}
\Li(((\bw \sh x^{\sh n}) \sh x)y;1) = n! \sum_{\varepsilon_1+\cdots+\varepsilon_r=n} \binom{k_1+\varepsilon_1-1}{\varepsilon_1}\cdots\binom{k_r+\varepsilon_r-1}{\varepsilon_r}\zeta(k_1+\varepsilon_1,\ldots,k_r+\varepsilon_r)
\end{equation}
by induction on $n$.

In case of $n=0$, the equation is $\Li(\bw y;1)=\zeta(k_1,\ldots,k_r)$. This is definition of $\Li$.

In case of $n=1$, clearly
\begin{equation}
(\bw \sh x)y = \sum_{j=1}^r k_j x^{k_1-1}y\cdots x^{k_j}y \cdots x^{k_r-1}y, \label{case_of_n=1}
\end{equation}
then the equation satisfies.

In general case, we assume that the lemma is satisfied on $n$.

Because of the hypothesis of induction and (\ref{case_of_n=1}), we get
\begin{align*}
&\Li((\bw \sh x^{\sh (n+1)})y;1) = \Li(((\bw \sh x^{\sh n}) \sh x)y;1) \\
& \quad = n! \sum_{\varepsilon_1+\cdots+\varepsilon_r=n}\binom{k_1+\varepsilon_1-1}{\varepsilon_1}\cdots\binom{k_r+\varepsilon_r-1}{\varepsilon_r} \sum_{j=1}^r (k_j+\varepsilon_j)\zeta(k_1+\varepsilon_1,\ldots,k_j+\varepsilon_j+1,\cdots,k_r+\varepsilon_r) \\
& \quad = n! \sum_{\varepsilon_1+\cdots+\varepsilon_r=n} \sum_{j=1}^r \binom{k_1+\varepsilon_1-1}{\varepsilon_1} \cdots (\varepsilon_j+1)\binom{k_j+\varepsilon_j}{\varepsilon_j+1} \cdots \binom{k_r+\varepsilon_r-1}{\varepsilon_r} \zeta(k_1+\varepsilon_1,\ldots,k_j+\varepsilon_j+1,\ldots,k_r+\varepsilon_r) \\
& \quad = n! \sum_{\varepsilon_1+\cdots+\varepsilon_r=n+1} \binom{k_j+\varepsilon_j-1}{\varepsilon_j} \cdots \binom{k_r+\varepsilon_r-1}{\varepsilon_r} \left(\sum_{j=1}^r \varepsilon_j\right)\zeta(k_1+\varepsilon_1,\ldots,k_r+\varepsilon_r)\\
& \quad = (n+1)! \sum_{\varepsilon_1+\cdots+\varepsilon_r=n+1} \binom{k_j+\varepsilon_j-1}{\varepsilon_j} \cdots \binom{k_r+\varepsilon_r-1}{\varepsilon_r} \zeta(k_1+\varepsilon_1,\ldots,k_r+\varepsilon_r)
\end{align*}

Therefore, the lemma is proved.

\end{proof}

\paragraph{Examples}

We compute some examples. First,
\begin{align}
&\frac{\ds \Gamma(1+(\lambda_2+\lambda_3))\Gamma(1-(\lambda_1+\lambda_2))}{\ds \Gamma(1+\lambda_3)\Gamma(1-\lambda_1)} = \exp(\sum_{n \ge 2} \frac{(-1)^n \zeta(n)}{n} ((\lambda_2+\lambda_3)^n-\lambda_3^n))\exp(\sum_{n \ge 2} \frac{\zeta(n)}{n} ((\lambda_1+\lambda_2)^n-\lambda_1^n)) \notag \\
&= \exp(\lambda_2 \sum_{n \ge 2}(-1)^n\zeta(n)\lambda_3^{n-1} + \lambda_2^2 \sum_{n \ge 2}\frac{(-1)^n (n-1)\zeta(n)}{2}\lambda_3^{n-2}+O(\lambda_2^3)) \notag \\
&\phantom{=} \times \exp(\lambda_2 \sum_{n \ge 2}\zeta(n)\lambda_1^{n-1} + \lambda_2^2 \sum_{n \ge 2}\frac{(n-1)\zeta(n)}{2}\lambda_1^{n-2}+O(\lambda_2^3)) \notag \\
&= 1 + \lambda_2\left(\sum_{n \ge 2} (-1)^n \zeta(n)\lambda_3^{n-1} + \sum_{n \ge 2} \zeta(n)\lambda_1^{n-1}\right) \notag \\
&\phantom{=} + \lambda_2^2 \left(\sum_{n \ge 2}\frac{(-1)^n (n-1)}{2}\zeta(n)\lambda_3^{n-2} + \sum_{n \ge 2}\frac{(n-1)}{2}\zeta(n)\lambda_1^{n-2} + \frac{1}{2}(\sum_{n \ge 2}(-1)^n\zeta(n)\lambda_3^{n-1})^2\right. \notag \\
&\phantom{= + \lambda_2^2 \bigg(} + \left.\frac{1}{2}(\sum_{n \ge 2}\zeta(n)\lambda_1^{n-1})^2 + \sum_{m \ge 2}\sum_{n \ge 2}(-1)^n \zeta(n)\zeta(m)\lambda_1^{m-1}\lambda_3^{n-1} \right) + O(\lambda_2^3) .
\end{align}

\begin{enumerate}

\item {\bf Case of $m=0$}

In this case, constant term of the formula (\ref{connection_12}) with respect to $\lambda^2$ is trivial: $0 = 0$.

\item {\bf Case of $m=1$}

If $l=0$ and $n \ge 1$, the coefficients of $\lambda_2\lambda_3^n$ are both $(-1)^{n+1}\zeta(n+1)$.

If $n=0$ and $l \ge 1$, we get easiest case of duality formula
\begin{equation*}
\zeta(2,\underbrace{1,\ldots,1}_{l-1 \text{ times}}) = \zeta(l+1)
\end{equation*}
as coefficients of $\lambda_1^{l}\lambda_2$.

If $l,m \ge 1$, formula \eqref{connection_12} degenerates as
\begin{equation}
\zeta(n+2,\underbrace{1,\ldots,1}_{l-1 \text{ times}}) + \sum_{j=0}^{n-1}\Li(x^{n-j}y^lx^{j+1};1)+\Li(x^{n-j}y^{l+1}x^j;1).
\end{equation}
By Lemma \ref{lem:culculus_hH},
\begin{align*}
\sum_{j=0}^{n-1}&\Li(x^{n-j}y^lx^{j+1};1) = \sum_{j=0}^{n-1}(-1)^{j+1}\sum_{\varepsilon=0}^{j+1}\binom{n-j+\varepsilon}{\varepsilon}\sum_{\varepsilon_2+\cdots+\varepsilon_l=j+1-\varepsilon}\zeta(n-j+1+\varepsilon,1+\varepsilon_2,\cdots,1+\varepsilon_l) \\
&= \sum_{j=0}^{n}(-1)^{j}\sum_{\varepsilon=0}^{j}\binom{n-j+\varepsilon+1}{\varepsilon}\sum_{\varepsilon_2+\cdots+\varepsilon_l=j-\varepsilon}\zeta(n-j+2+\varepsilon,1+\varepsilon_2,\cdots,1+\varepsilon_l) - \zeta(n+2,1,\ldots,1)\\
&= \sum_{j=0}^{n}(-1)^{j}\sum_{\varepsilon=0}^{j}\binom{n+1-\varepsilon}{j-\varepsilon}\sum_{\varepsilon_2+\cdots+\varepsilon_l=\varepsilon}\zeta(n+2-\varepsilon,1+\varepsilon_2,\cdots,1+\varepsilon_l) - \zeta(n+2,1,\ldots,1)\\
&= \sum_{\varepsilon=0}^{n}\sum_{\varepsilon_2+\cdots+\varepsilon_l=\varepsilon} \left(\sum_{j=0}^{n}(-1)^{j}\binom{n+1-\varepsilon}{j-\varepsilon}\right)\zeta(n+2-\varepsilon,1+\varepsilon_2,\cdots,1+\varepsilon_l) - \zeta(n+2,1,\ldots,1)\\
&= (-1)^n \sum_{\substack{k_1+\cdots+k_l=n+l+1\\k_1 \ge 2}}\zeta(k_1,\ldots,k_l)- \zeta(n+2,1,\ldots,1)
\end{align*}
and by the same way,
\begin{equation*}
\sum_{j=0}^{n-1}\Li(x^{n-j}y^{l+1}x^{j};1) = -(-1)^n \sum_{\substack{k_1+\cdots+k_{l+1}=n+l+1\\k_1 \ge 2}}\zeta(k_1,\ldots,k_{l+1}).
\end{equation*}

Then we get
\begin{align}
\sum_{\substack{k_1+\cdots+k_l=n+l+1\\k_1 \ge 2}}\zeta(k_1,\ldots,k_l)-\sum_{\substack{k_1+\cdots+k_{l+1}=n+l+1\\k_1 \ge 2}}\zeta(k_1,\ldots,k_{l+1}) = 0.
\end{align}

This show that $\forall n,l \ge 1$,
\begin{equation}
\sum_{\substack{k_1+\cdots+k_l=n+l+1\\k_1 \ge 2}}\zeta(k_1,\ldots,k_l) = \zeta(n+l+1).
\end{equation}

This is the sum formula of multiple zeta values.

\item {\bf Case of $m=2, n=0$}

By similar computing, we get
\begin{equation}
(l+1)\zeta(2,\underbrace{1,\ldots,1}_{l \text{ times}}) - \zeta(3,\underbrace{1,\ldots,1}_{l-1 \text{ times}})= \frac{l+1}{2}\zeta(l+2)+\frac{1}{2}\sum_{\substack{i,j \ge 0\\i+j=l-2}}\zeta(i+2)\zeta(j+2)
\end{equation}

This formula is equal to one of the Euler's formula(cf. Zagier\cite{Z1}) up to duality.

\item {\bf Case of $m=2, n=1$}

\begin{align}
&2G_0(l+3,l,1)+G_0(l+3,l,2)-(l+2)G_0(l+3,l+1,1)-lG_0(l+3,l+1,2)+(l+1)G_0(l+3,l+2,1) \notag \\
&= \begin{cases}-\zeta(3)&\quad (l=0)\\ \zeta(2)\zeta(l+1) &\quad (l \ge 1)\end{cases}
\end{align}

\end{enumerate}

In general case, computing is very complicated. We can't write all coefficients expricitly, but we conjecture that all coefficients of formula \eqref{connection_12}'s LHS are $\Z$ linear combinations of $G_0(k,n,s;1)$ and the formula \eqref{connection_12} is also equal to Ohno-Zagier ralation by change of values.

\paragraph{}

We don't compute a solution of the hypergeometric equation singular at $z=1$ yet, but this can be expressed in similar way and limit of $(2,1)$ component of connection relation(\ref{connection01}) as $z \to 0$ is converge to "dual" of results of section \ref{subsec:OtherSolution:singular}.

However, limit of $(1,2)$ component as $z \to 1$, limit of $(2,1)$ component as $z \to 0$ and limit of $(2,2)$ component as $z \to 0,1$ are apparently diverge. Computing these limits are very difficult and needed careful analysis. We don't treat these cases in this paper.

\subsection{Solutions of hypergeometric equation in neighborhood of $\infty$} \label{subsec:OtherSolution:infty}

In the last of this paper, we compute $\Phi^{-1}_{\infty}$: inverse of fundamental solution matrix in neighborhood of $\infty$.

Let $u=1/z$. The equation \eqref{HGeq2} can be rewrite as:

\begin{align}
\frac{d}{du} {}^t \Phi_{\infty}^{-1}(u) &= (\lambda_1 \theta''_1 + \lambda_2 \theta''_2 + \lambda_3 \theta''_3) {}^t \Phi_{\infty}^{-1}(u)\label{Phi_infty_inv_eq}
\end{align}
\begin{align*}
\theta''_1 &= \frac{1}{u}\begin{pmatrix}0 & 0 \\ 1 & -1\end{pmatrix} + \frac{1}{1-u}\begin{pmatrix}0 & 0 \\ 0 & -1\end{pmatrix} \\
\theta''_2 &= \frac{1}{1-u}\begin{pmatrix}0 & 0 \\ 0 & -1\end{pmatrix} \\
\theta''_3 &= \frac{1}{u}\begin{pmatrix}0 & 1 \\ 0 & 1\end{pmatrix} + \frac{1}{1-u}\begin{pmatrix}0 & 1 \\ 0 & 0\end{pmatrix}
\end{align*}

Asymptotic condition of $\Phi_{\infty}^{-1}$ is
\begin{equation}
\Phi_{\infty}^{-1}(u) \to \frac{\beta}{\beta-\alpha} \begin{pmatrix}u^{-\alpha}&-\frac{\alpha}{\beta}u^{-\beta}\\ u^{-\alpha}&-u^{-\beta}\end{pmatrix} \quad (u \to 0)
\end{equation}
by direct computing.

We denote
\begin{equation}
\Phi_{\infty}^{-1}(u)=\begin{pmatrix}\psi^{(\infty)}_{11}(u) & \psi^{(\infty)}_{12}(u) \\ \psi^{(\infty)}_{21}(u) & \psi^{(\infty)}_{22}(u)\end{pmatrix}.
\end{equation}

In the same way as above, we can construct the solutions of \eqref{Phi_infty_inv_eq} which are made in iterated integral starting from $\vect{1}{0}$ and $\vect{1}{0}$.

\begin{lem}

Define the transformation $T_{\infty}:\{\text{finite sequences of}\{1,2,3\}\} \to \fH$ as:
\begin{enumerate}
\item $T_{\infty}(\emptyset)=1$
\item $T_{\infty}(1,3,\bmu)=(xy-yx)T_1(\bmu)$
\item $T_{\infty}(1,\bmu)=-(x+y)T_1(\bmu) \qquad$ if $\bmu$ doesn't start with 3
\item $T_{\infty}(2,\bmu)=-yT_1(\bmu)$
\item $T_{\infty}(3,\bmu)=xT_1(\bmu)$.
\end{enumerate}

Then the solution of \eqref{Phi_infty_inv_eq} which is made in iterated integral starting from $\vect{1}{0}$ is
\begin{equation}
\vect{1}{0} +\vect{\displaystyle \sum_{r=0}^{\infty} \sum_{\mu_1,\ldots,\mu_r = 1,2,3} \lambda_3 \lambda_{\mu_1}\cdots\lambda_{\mu_r} \lambda_1 \Li((x+y)T_{\infty}(\mu_1,\ldots,\mu_r)x;u)}{\displaystyle \sum_{r=0}^{\infty} \sum_{\mu_1,\ldots,\mu_r = 1,2,3} \lambda_{\mu_1}\cdots\lambda_{\mu_r} \lambda_1 \Li(T_{\infty}(\mu_1,\ldots,\mu_r)x;u)} \label{solution_infty_10}
\end{equation}
and starting from $\vect{0}{1}$ is
\begin{equation}
\vect{0}{1} +\vect{\displaystyle \sum_{r=0}^{\infty} \sum_{\mu_1,\ldots,\mu_r = 1,2,3} \lambda_3 \lambda_{\mu_1}\cdots\lambda_{\mu_r} \Li((x+y)T_{\infty}(\mu_1,\ldots,\mu_r);u)}{\displaystyle \sum_{r=1}^{\infty} \sum_{\mu_1,\ldots,\mu_r = 1,2,3} \lambda_{\mu_1}\cdots\lambda_{\mu_r} \Li(T_{\infty}(\mu_1,\ldots,\mu_r);u)} \label{solution_infty_01}
\end{equation}

\end{lem}

By the lemma and asymptotic condition, the solution $\vect{\psi^{(\infty)}_{11}(u)}{\psi^{(\infty)}_{12}(u)}$ can be constructed as $\frac{\lambda_1}{\lambda_1+\lambda_3}$(\eqref{solution_infty_10}$+$\eqref{solution_infty_01}) and $\vect{\psi^{(\infty)}_{21}(u)}{\psi^{(\infty)}_{22}(u)}$ to $\frac{1}{\lambda_1+\lambda_3}$($\lambda_3$ \eqref{solution_infty_10}$+\lambda_1$ \eqref{solution_infty_01}). Then we get the next proposition.

\begin{prop}

\begin{equation}
\vect{\psi^{(\infty)}_{11}(u)}{\psi^{(\infty)}_{12}(u)} = \vect{\ds \sum_{l,m,n}\left( \sum_{\bmu \in J(l-1,m,n-1)} \Li((x+y)T_{\infty}(\bmu)x;u)+\sum_{\bmu \in J(l,m,n-1)} \Li((x+y)T_{\infty}(\bmu);u)\right)\lambda_1^l\lambda_2^m\lambda_3^n+1}{\ds \sum_{l,m,n}\left( \sum_{\bmu \in J(l-1,m,n)} \Li(T_{\infty}(\bmu)x;u)+\sum_{\bmu \in J(l,m,n)} \Li(T_{\infty}(\bmu);u)\right)\lambda_1^l\lambda_2^m\lambda_3^n}
\end{equation}

\end{prop}

Comparing these results and connection relation \eqref{connection0infty}, we get some relations. But in general the formulation is very complicated. The connection relation of $0$ and $\infty$ is also regarded as Euler's inversion formula. Indeed computing $(1,1)$ component of \eqref{connection0infty} and comparing coefficient of $\lambda_3^n$, we get the next formula.

\begin{prop}
\begin{align}
&\frac{\lambda_1}{\lambda_1+\lambda_3}\left(\frac{\log^n \frac{1}{z}}{n!} - \sum_{k=0}^{n-1}(\Li_{n-k}(z)+(-1)^{n-k}\Li_{n-k}(\frac{1}{z}))\frac{\log^{k} \frac{1}{z}}{k!}\right) \notag\\
&= \frac{\lambda_1}{\lambda_1+\lambda_3}\exp(\pi i \lambda_3)\frac{\pi \lambda_3}{\sin \pi \lambda_3} = \frac{\lambda_1}{\lambda_1+\lambda_3} \sum_{n=0}^{\infty} B_n \frac{(2 \pi i \lambda_3)^n}{n!}
\end{align}
Now $B_n$ are Bernoulli numbers.

Especially if $n:$even, considering $z \to 1$, we get
\begin{equation*}
-2\zeta(n)=B_n\frac{(2 \pi i)^n}{n!}.
\end{equation*}

This is Euler's famous formula.

\end{prop}

\paragraph{}

\begin{flushleft}
\hspace{3cm}Shu OI\\
\hspace{3cm}Department of Mathematical Sciences\\
\hspace{3cm}School of Science and Engineering\\
\hspace{3cm}Waseda University\\
\hspace{3cm}Okubo 3-4-1, Shinjuku-ku, Tokyo 169-8555, Japan\\
\hspace{3cm}{\bf e-mail:} {\tt shu@gm.math.waseda.ac.jp}\\
\end{flushleft}

\end{document}